\documentclass[11pt,reqno]{amsart}
\usepackage[T1]{fontenc}
\usepackage{lmodern}
\usepackage[a4paper,margin=29mm]{geometry}
\usepackage{amsmath,amssymb,mathtools,mathrsfs,tikz-cd}
\usepackage{microtype,booktabs,longtable,array,enumitem,xcolor,xurl}
\definecolor{linkblue}{RGB}{30,65,100}
\usepackage[colorlinks=true,linkcolor=linkblue,citecolor=linkblue,urlcolor=linkblue]{hyperref}
\usepackage{bookmark}
\newcommand{\kk}{\Bbbk}\newcommand{\ZZ}{\mathbb Z}\newcommand{\NN}{\mathbb N}
\newcommand{\OO}{\mathcal O}\newcommand{\CC}{\mathscr C}\newcommand{\RR}{\mathcal R}\newcommand{\DD}{\mathcal D}
\newcommand{\Hom}{\operatorname{Hom}}\newcommand{\End}{\operatorname{End}}\newcommand{\Ext}{\operatorname{Ext}}
\newcommand{\Dist}{\operatorname{Dist}}\newcommand{\Ind}{\operatorname{Ind}}\newcommand{\Perf}{\operatorname{Perf}}
\newcommand{\add}{\operatorname{add}}\newcommand{\rad}{\operatorname{rad}}
\newcommand{\fl}{\mathrm{fl}}\newcommand{\pf}{\mathrm{pf}}\newcommand{\op}{\mathrm{op}}
\newcommand{\BJ}{\mathsf B_J}\newcommand{\BK}{\mathsf B_K}\newcommand{\TJ}{\mathbb T_J}\newcommand{\GK}{\mathsf G_K}\newcommand{\GJ}{\mathsf G_J}
\newtheorem{theorem}{Theorem}[section]
\newtheorem{proposition}[theorem]{Proposition}
\newtheorem{lemma}[theorem]{Lemma}
\newtheorem{corollary}[theorem]{Corollary}
\theoremstyle{definition}

\newtheorem{example}[theorem]{Example}
\theoremstyle{remark}

\numberwithin{equation}{section}
\setlist[enumerate,1]{label=\textup{(\roman*)},leftmargin=*,itemsep=3pt}
\allowdisplaybreaks[2]
\makeatletter
\renewcommand{\@maketitle}{%
  \normalfont\normalsize
  \@adminfootnotes
  {\centering\baselineskip14pt\bfseries
    \uppercasenonmath\@title\@title\par}%
  \vskip12pt
  {\centering\footnotesize\MakeUppercase{\authors}\par}%
  \vskip14pt
  \box\abstractbox
  \vskip10pt
}
\newcommand{\makefrontpage}{%
  \thispagestyle{firstpage}%
  \vspace*{\fill}%
  \maketitle
  \begingroup
      \centering
        \begin{minipage}{.94\textwidth}
          \normalsize
          \linespacing=\baselineskip
          \def\l@section{\@tocline{1}{0pt}{0pt}{}{}}%
          \def\l@subsection{\@tocline{2}{0pt}{0pt}{5pc}{}}%
          \tableofcontents
          \par\removelastskip
        \end{minipage}%
      \par
  \endgroup
  \vspace*{\fill}%
  \markboth{DERIVED RINGEL DUALITY}{FINITE SUBQUOTIENTS OF MODULAR CATEGORY $\OO$}%
  \clearpage
}
\makeatother

\hypersetup{pdftitle={Derived Ringel Duality through Finite Subquotients of Modular Category O},pdfauthor={Chenhan Zhu},pdfsubject={Modular category O, finite Serre subquotients and derived Ringel duality}}
\title[Derived Ringel duality in modular category $\OO$]
{Derived Ringel duality through finite subquotients\\
of modular category $\OO$}
\author{Chenhan Zhu}
\address{School of Mathematical Sciences, East China Normal University, Shanghai 200241, China}
\email{chenhanzhumath@gmail.com}
\date{}
\begin{document}
\begin{abstract}
We study derived Ringel duality for modular category $\OO$ on
downward-closed, quasi-bounded supports. We reconstruct global tilting
objects and their morphisms from compatible images in finite
highest-weight Serre subquotients, and realize the complete Ringel
algebra as an inverse limit of finite-dimensional algebras. The duality
identifies the bounded homotopy category of finitely generated
projectives over the associated locally unital algebra with the
bounded derived category of the Ringel dual category. We characterize tiltings
recovered, together with their outgoing morphisms, by a single finite
subquotient, and give a sharp criterion for equivalence under Frobenius
restriction. For $SL_2$, we compute tilting endomorphism algebras and
their quotient maps, and obtain sharp recovery bounds for finite
tilting families.
\end{abstract}
\makefrontpage
\addtocontents{toc}{\protect\setcounter{tocdepth}{1}}
% BEGIN introduction.tex
\section{Introduction}\label{sec:introduction}

Modular category $\OO$ provides a highest-weight setting for
representations of the distribution algebra of a semisimple algebraic
group in positive characteristic. Its ordinary Verma modules are
generally infinite dimensional, and the tilting objects considered here
may have infinite Verma flags. Finite highest-weight Serre subquotients
make parts of this structure accessible through finite-dimensional
algebras. Our aim is to understand how these finite categories determine
global tilting objects, their morphisms and the associated Ringel
duality.

Fiebig constructs subquotients associated with locally closed sets of
weights and proves equivalences under suitable weight shifts
\cite[Section~17 and Theorem~25]{Fiebig}. Andersen develops the theory
of tilting modules with finite Verma flags \cite[Section~7]{Andersen}.
Brundan--Stroppel's semi-infinite Ringel duality places tiltings with
possibly infinite flags in the framework of upper and lower finite
highest-weight categories \cite[Theorem~4.27]{BS}. We use these results
to describe Ringel duality through the quotient functors between
Fiebig's finite subquotients.

We fix a downward-closed support in which only finitely many weights
lie above any given weight. Compatible tilting objects in the finite
subquotients reconstruct global tiltings, and compatible morphisms
reconstruct their Hom spaces. We determine the kernels of the induced
maps between tilting endomorphism algebras and identify their inverse
limit with the complete Ringel algebra. These constructions yield an
equivalence between the bounded homotopy category of finitely generated
projectives in the supported category and the bounded derived category
of its Ringel dual. This realizes the homotopy-category formulation
described in \cite[Section~2.4, equation~(2.3) and the following
display]{FG}, compatibly with the finite subquotients.

This compatibility gives a finite setting for homological calculations.
Any finite collection of objects, morphisms and exact triangles in the
bounded derived category of the Ringel dual comes from the bounded
derived category of a finite highest-weight subcategory. Extension
groups between objects of that subcategory and compositions of derived
morphisms can therefore be computed there.
We also determine exactly which global tiltings have all their
outgoing morphisms recovered by one finite subquotient. These are
precisely the tiltings with finite ordinary Verma flags, and precisely
those represented by bounded complexes of finitely generated
projectives.

Frobenius restriction gives a further finite realization. We establish
a sharp bound under which restriction to a Frobenius kernel, with the
rational torus retained, induces an equivalence of the finite
subquotient categories. For $SL_2$, the construction gives explicit
tilting endomorphism algebras, their quotient maps, and sharp bounds
for preserving all morphisms in finite tilting families. These
calculations illustrate how the categorical description determines
both algebra structures and the amount of finite data needed to
recover them.

\subsection{Finite subquotients and Frobenius restriction}
Let $\kk$ be algebraically closed of characteristic $p>0$, and let $G$
be split, semisimple and simply connected, with its fixed split
$\mathbb F_p$-form. Fix a maximal torus
$T$ and positive Borel. Write $U=\Dist(G)$, $X$ for the
weight lattice and $\Pi$ for the simple roots, with
$\mu\leq\nu$ when $\nu-\mu\in\NN\Pi$. The rational torus is retained
throughout.

Give $X$ Fiebig's topology, whose open sets are the downward-closed
subsets for this order \cite[Section~10, Definition~3]{Fiebig}.
Fix an open support $J\subseteq X$ such that
$\{\nu\in J:\mu\leq\nu\}$ is finite for every $\mu$. Let $\RR_J$
be the category of integral $U$-weight modules supported in $J$ with
finite-dimensional weight spaces. A \emph{window} is a nonempty
finite subset $K$ closed in $J$. Write $\CC_K$ for the finite-length
part of $\OO^J/\OO^{J\setminus K}$, where $\OO^J$ allows arbitrary
weight-space dimensions. The exact quotient functor
$q_K:\RR_J\to\CC_K$ has kernel the modules supported in
$J\setminus K$. Restriction to $G_rT$ induces
$R_{r,K}:\CC_K\to\CC_{r,K}$.

The weight projective $Q^J(\mu)$ represents $M\mapsto M_\mu$.
Let $A_J$ be their locally unital endomorphism algebra with opposite
composition, and put
$e_K=\sum_{\mu\in K}1_\mu$, $A_K=e_KA_Je_K$. Define
\[
 d(K)=\max\bigl(\{0\}\cup
 \{n_i:\nu-\mu=\textstyle\sum_i n_i\alpha_i,\ 
                         \mu,\nu\in K,\ \mu\leq\nu\}\bigr).
\]

\begin{theorem}[Finite-window realization]\label{thm:intro-window}
The weight-projective equivalences
\[
 \RR_J\simeq A_J\text{-}\mathrm{mod}_{\mathrm{lfd}},
 \qquad \CC_K\simeq A_K\text{-}\mathrm{mod}_{\mathrm{fd}}
\]
identify $q_K$ with $e_K(-)$, and
\[
 A_K\cong
 \End_{\CC_K}\Bigl(\bigoplus_{\mu\in K}q_KQ^J(\mu)\Bigr)^{\op}.
\]
For every $r\geq1$,
\[
 R_{r,K}\text{ is an equivalence}
 \quad\Longleftrightarrow\quad p^r>d(K).
\]
These realizations commute with the quotient functors
$q_K^L:\CC_L\to\CC_K$ for $K\subseteq L$ and restriction
between Frobenius kernels.
\end{theorem}

Here locally finite dimensional means finite dimensional at each weight
idempotent. Theorems~\ref{thm:raw} and~\ref{thm:frobenius} prove the
statement using PBW bases. For the canonical support
$J=\mathord\downarrow K$, dominant $\gamma$ and $p^\ell>d(K)$,
Section~\ref{sec:frobenius} compares Fiebig's shift
$\tau=p^\ell\gamma$ with tensoring by the torus character $\kk_\tau$
at Frobenius levels $r\leq\ell$.

\subsection{Global tiltings and finite recovery}
Write $\Delta(\lambda)$ and $\nabla(\lambda)$ for the ordinary Verma
module and its restricted dual, and put
$\Delta_K(\lambda)=q_K\Delta(\lambda)$ and
$\nabla_K(\lambda)=q_K\nabla(\lambda)$ for $\lambda\in K$.
The upper finite highest-weight
category $\RR_J$ has a self-dual indecomposable tilting
$T(\lambda)$ \cite[Theorem~4.18]{BS}, whose ordinary flags may
be infinite. Its image in a window retaining $\lambda$ is
$T_K(\lambda)$. The multiplicities
\[
 t_\lambda(\mu)=(T_K(\lambda):\Delta_K(\mu)),\qquad
 S(T(\lambda))=\{\mu:t_\lambda(\mu)>0\}
\]
are independent of the window containing both labels. We call
$S(T(\lambda))$ its standard support.

For a window $K$ containing $\lambda$, let $j_{K,!}$ and $j_{K,*}$
denote the left and right adjoints of $q_K$, respectively, following
\cite[Section~2.5]{BS}. They give a canonical map
\[
 c_K:j_{K,!}T_K(\lambda)\longrightarrow j_{K,*}T_K(\lambda)
\]
which becomes the identity under $q_K$.

\begin{theorem}[Finite recovery]\label{thm:intro-recovery}
For $T=T(\lambda)$, the following are equivalent:
\begin{enumerate}
\item $T$ is finitely generated, equivalently finitely presented.
\item $T$ has finite ordinary standard and costandard flags.
\item $S(T)$ is finite.
\item $\End_U(T)$ is finite dimensional.
\item Some window $K$ recovers every outgoing morphism:
\[
 \Hom_U(T,N)\xrightarrow{\sim}
       \Hom_{\CC_K}(q_KT,q_KN)\qquad(N\in\RR_J).
\]
\end{enumerate}
For a fixed window retaining $\lambda$, the map $c_K$ is an
isomorphism, a monomorphism or an epimorphism precisely when
$S(T)\subseteq K$. If $S(T)$ is finite, its upper closure in $J$
is the least recovering window.

More generally, for a finite set $\Gamma\subseteq K$, the functor $q_K$ is
fully faithful on $\add(\bigoplus_{\lambda\in\Gamma}T(\lambda))$
if and only if $\bigcup_{\lambda\in\Gamma}S(T(\lambda))\subseteq K$.
\end{theorem}

Theorems~\ref{thm:recovery}--\ref{thm:adjoint-comparison} and
Corollary~\ref{cor:family-quotient} give these criteria. The underlying
mixed Hom calculation is
\[
 \Hom_U(T(\lambda),T(\mu))
 \cong\varprojlim_K\Hom_{\CC_K}(T_K(\lambda),T_K(\mu)).
\]
Here the limit runs over windows containing both $\lambda$ and $\mu$,
with transition maps induced by $q_K^L$.
Each quotient-induced map $f\mapsto q_K(f)$ is surjective and is
an isomorphism exactly when
the common standard support
$\{\nu:t_\lambda(\nu)t_\mu(\nu)>0\}$ lies in $K$
(Proposition~\ref{prop:tilting-hom}).

\subsection{The complete Ringel algebra and derived duality}
Set
\[
 \TJ=\bigoplus_{\lambda\in J}T(\lambda),\qquad
 \Theta_K=q_K\TJ,\qquad
 \BK=\End_{\CC_K}(\Theta_K)^{\op},\qquad
 \BJ=\End_U(\TJ)^{\op}.
\]
Quasi-boundedness makes $\TJ$ pointwise finite.
For $\lambda\in K$, $\varepsilon_{\lambda,K}$ is the endomorphism
of $\Theta_K$ equal to the identity on $T_K(\lambda)$ and zero on
the other summands, viewed in $\BK$. The algebra maps
$\pi_K:\BJ\to\BK$, $f^{\op}\mapsto(q_Kf)^{\op}$, define the
profinite topology on $\BJ$: the ideals $\ker\pi_K$ form a
neighbourhood basis at zero, with each $\BK$ discrete.
Section~\ref{sec:ringel} proves completeness and identifies this
topology with that used for semi-infinite Ringel duality.

\begin{theorem}[Ringel duality in finite windows]\label{thm:intro-ringel}
For $K\subseteq L$, the functor $q_K^L$ induces a surjective algebra map
$\pi_{L,K}:\mathsf B_L\to\BK$, $f^{\op}\mapsto(q_K^Lf)^{\op}$, with
\[
 \ker\pi_{L,K}=\mathsf B_L
       \Bigl(\sum_{\nu\in L\setminus K}\varepsilon_{\nu,L}\Bigr)
       \mathsf B_L.
\]
There are topological algebra isomorphisms
\[
 \BJ\cong\varprojlim_K\BK
 \cong\varprojlim_{\substack{(K,r)\\p^r>d(K)}}
     \End_{\CC_{r,K}}(R_{r,K}\Theta_K)^{\op}.
\]
Let $\DD_J$ consist of finite-dimensional $\BJ$-modules whose
action factors through some $\BK$. It is lower finite highest weight,
with label order $J^{\op}$. Via inflation along $\BJ\to\BK$, the
categories $\DD_K=\BK\text{-}\mathrm{mod}_{\mathrm{fd}}$ identify
with the Serre subcategories on simple labels $K$. They exhaust
$\DD_J$. Writing $L^\vee,\Delta^\vee,\nabla^\vee$ for the simple,
standard and costandard objects of this dual category, respectively,
we have
\[
 t_\lambda(\mu)
   =[\Delta^\vee(\mu):L^\vee(\lambda)]
   =[\nabla^\vee(\mu):L^\vee(\lambda)].
\]
\end{theorem}

The categorical duality is Brundan--Stroppel's \cite{BS};
Theorem~\ref{thm:ringel}, Corollary~\ref{thm:dual-reciprocity}
and Corollary~\ref{cor:frobenius-ringel} identify its maps in the
window realization. The covariant Ringel functor $\GJ$ is the
continuous dual of $\Hom_U(-,\TJ)$, with finite counterpart
$\GK=\Hom_{\CC_K}(-,\Theta_K)^*$. Its images
$\mathsf W_K=\GK(\bigoplus_{\mu\in K}q_KQ^J(\mu))$
recover $A_K$:
\[
 A_K\cong\End_{\DD_K}(\mathsf W_K)^{\op}
 \qquad\text{(Proposition~\ref{prop:raw-bidual}).}
\]

Write $i_K:\DD_K\hookrightarrow\DD_J$ for inflation: a $\BK$-module
keeps its vector space, with $\BJ$ acting through the quotient map.
Let $\Perf(A)=K^b(\operatorname{proj}_{\mathrm{fg}}A)$ be the bounded
homotopy category of finitely generated projectives.
Write $D^b$ for the bounded derived category and $\mathbf L$ for
a total left-derived functor. On $\Perf(A_J)$, $\mathbf G_J$
applies $\GJ$ termwise. The exact functor $i_K$ induces $D^b(i_K)$
termwise as well; Section~\ref{sec:derived} fixes the conventions.

\begin{theorem}[Derived realization and heredity layers]\label{thm:intro-derived}
There is a naturally commutative square
\[
\begin{tikzcd}[column sep=large,row sep=large]
 \Perf(A_K)\arrow[r,"\mathbf L\GK","\sim"']
             \arrow[d,"\mathbf L j_{K,!}"']&
 D^b(\DD_K)\arrow[d,"D^b(i_K)"]\\
 \Perf(A_J)\arrow[r,"\mathbf G_J"',"\sim"]&
 D^b(\DD_J).
\end{tikzcd}
\]
Both vertical functors are fully faithful. As $K$ varies, their
essential images exhaust the global categories. Each $\DD_K$ has a heredity
filtration with derived strata $D^b(\kk\text{-}\mathrm{mod}_{\mathrm{fd}})$.
For $T(\lambda)$, the conditions of Theorem~\ref{thm:intro-recovery}
are also equivalent to $T(\lambda)$ being perfect and to its dual
injective hull $I^\vee(\lambda)$ having finite length.
\end{theorem}

Theorem~\ref{thm:derived}, Corollary~\ref{cor:perfect-recovery} and
Proposition~\ref{prop:heredity-layer} prove the assertions through
finite derived Ringel duality and heredity ideals. 

\subsection{Rank-one towers and mixed algebras}
For $G=SL_2$, identify the positive root with $2$ and set
$K_{\lambda,N}=\{\lambda-2n:0\leq n\leq N\}$.
Expand $\lambda+1=\sum_{i\geq0}a_ip^i$ in $\ZZ_p$, and put
$I_\lambda=\{i:a_i\ne0\}$. For finite $A\subseteq I_\lambda$, write
$d_A=\sum_{i\in A}a_ip^i$ and $x_A=\prod_{i\in A}x_i$.

\begin{theorem}[The depth-indexed rank-one tower]\label{thm:intro-rank-one}
For every prime $p$ and integer $\lambda$,
\[
 t_\lambda(\lambda-2n)=
 \begin{cases}
 1,&n=d_A\text{ for some finite }A\subseteq I_\lambda,\\
 0,&\text{otherwise}.
 \end{cases}
\]
Moreover,
\[
 \End_{\CC_{K_{\lambda,N}}}(T_{K_{\lambda,N}}(\lambda))
 \cong
 \kk[x_i:i\in I_\lambda]/(x_i^2,\ x_A:d_A>N).
\]
For $N\leq M$, the algebra map induced by $q_{K_{\lambda,N}}^{K_{\lambda,M}}$
keeps $x_A$ when $d_A\leq N$ and sends it to zero otherwise.
The inverse limit is $\End_U(T(\lambda))$.
The dimension profile determines $\lambda$; prefixes through depth
$p^m-1$ are isomorphic precisely when the labels agree modulo $p^m$.
The ordinary flags of $T(\lambda)$ are finite exactly when
$\lambda\geq-1$, in which case the least recovering depth is
$N=\lambda+1$.
\end{theorem}

The proof combines Andersen's tensor construction \cite{Andersen},
rational tilting morphisms \cite{TW}, and Fiebig's periodicity
\cite{Fiebig}; see Theorems~\ref{thm:rank-one-tower}
and~\ref{thm:depth-rigidity}. Numerical depth is essential:
all labels below $-1$ have isomorphic unfiltered topological
endomorphism algebras.

For a general interval $K=\{b,b-2,\ldots,a\}$, write
$t^K_\lambda(\eta)=(T_K(\lambda):\Delta_K(\eta))$, extended by zero
outside $K$, and set
$h_{\lambda,r}(\eta)=\sum_{j\geq0}t^K_\lambda(\eta+2jp^r)$.

\begin{theorem}[Mixed Frobenius recovery]\label{thm:intro-mixed}
For every $r\geq1$, restriction gives
\[
 \dim\Hom_{\CC_{r,K}}(R_{r,K}T_K(\lambda),R_{r,K}T_K(\mu))
       =\sum_{\eta\in K}h_{\lambda,r}(\eta)h_{\mu,r}(\eta).
\]
For a nonempty $\Gamma\subseteq K$, $R_{r,K}$ is fully faithful on
$\add(\bigoplus_{\lambda\in\Gamma}T_K(\lambda))$ precisely when
\[
 p^r>\max_{\lambda\in\Gamma}\frac{\lambda-a}{2}.
\]
The same condition is detected by the diagonal endomorphism rings.
For a global family, recovery through $R_{r,K}q_K$
holds precisely when this inequality and
$\bigcup_{\lambda\in\Gamma}S(T(\lambda))\subseteq K$ both hold.
\end{theorem}

Theorem~\ref{thm:mixed-gram} and
Corollaries~\ref{cor:sharp-family}--\ref{cor:composite-family} follow
by refining both ordinary flags into baby-standard flags. The two
flags are crucial: Example~\ref{ex:mixed} gives ordinary standards
whose diagonal rings are recovered while an extra mixed map appears.

Mixed products admit finite quiver presentations as well. Put
$N=(b-a)/2$, choose $\ell,u\geq1$ with
$H=p^\ell>\max(N,-a-1)$ and $Q=p^u>b+2H+1$,
and set $\tau=2H$, $c=Q+b+\tau$, $\alpha=Q+a+\tau$. Let
\[
 Z_c=\End_{SL_2}\Bigl(
       \bigoplus_{\substack{0\leq m\leq c\\m\equiv c\ (2)}}
                      T^{SL_2}(m)\Bigr),\qquad
 e_{<\alpha}=\sum_{\substack{0\leq m<\alpha\\m\equiv c\ (2)}}e_m,
\]
where the $T^{SL_2}(m)$ are rational tiltings and $e_m$ are their
summand idempotents.

\begin{theorem}[Finite mixed presentations]\label{thm:intro-presentation}
There is an algebra isomorphism
\[
 \BK\cong (Z_c/Z_ce_{<\alpha}Z_c)^{\op},
\]
where the vertex labelled by the rational highest weight
$m=Q+\lambda+\tau$ corresponds to $T_K(\lambda)$.
On every finite chain of nested intervals with the same upper bound,
common parameters make the transition maps the algebra quotient maps by
the additional lower-vertex ideals.
\end{theorem}

Theorem~\ref{thm:quiver-transport} transports the rational quiver
presentation of \cite{TW}. Example~\ref{ex:mixed-four} gives its
four-weight algebra, products and Frobenius recovery failure explicitly.

\subsection{Sources and organization}
Tilting multiplicities also
play a central role in the related work of Soergel \cite{Soergel}.
The heredity layers use \cite{Krause}, and the rank-one presentations
use \cite{TW,STWZ}.

Sections~\ref{sec:windows}--\ref{sec:ringel} develop recovery and
Ringel duality. Section~\ref{sec:examples} gives the rank-one towers,
mixed tilting presentations and sharp recovery bounds. A separate
notation companion collects the definitions and compares the
conventions used in the principal references.
% END introduction.tex

\addtocontents{toc}{\protect\setcounter{tocdepth}{2}}
% BEGIN foundations.tex
\section{Finite subquotients and Frobenius realization}\label{sec:windows}

Weight projectives describe the finite quotients and reconstruct
global objects from their compatible images.

\subsection{The supported category}
Let $\kk$ be an algebraically closed field of characteristic $p>0$, and
let $G$ be a split semisimple simply connected group. Fix a maximal torus
$T$ and a positive Borel subgroup. Write $X$ for the weight
lattice, $R^+$ for the positive roots and $\Pi=\{\alpha_i\}$ for the
simple roots. We use the root order
\[
 \mu\leq\lambda\quad\Longleftrightarrow\quad
 \lambda-\mu\in\NN\Pi.
\]
The distribution algebra $U=\Dist(G)$ has the triangular decomposition
$U=U^-U^0U^+$. Its positive Borel subalgebra is $B=U^0U^+$.
We use integral divided powers.

Following the modular category $\OO$ conventions of
\cite{Fiebig,Andersen}, we take $\OO$ to consist of integral weight
$U$-modules on which $U^+$ acts locally finitely. Each module carries
a rational torus action agreeing with the $U^0$-action and satisfying
$t(um)=(\operatorname{Ad}t)(u)(tm)$; morphisms preserve this structure.

Equip $X$ with Fiebig's topology, whose open sets are the
downward-closed subsets \cite[Section~10, Definition~3]{Fiebig}.
Fix an open subset $J\subseteq X$ such that
\[
 J_{\geq\mu}:=\{\nu\in J:\mu\leq\nu\}
 \quad\text{is finite for every }\mu\in X.
\]
This finiteness condition is called quasi-boundedness. Let $\OO^J$
be the full subcategory supported in $J$, and put
$\RR_J=\OO^J_{\pf}$, where the subscript requires each weight space
to be finite dimensional. A typical choice is
$J=\operatorname{down}(\lambda)$; all subquotients below use this fixed
support unless otherwise stated.

For $\lambda\in J$, the Verma and dual Verma modules are
\[
 \Delta(\lambda)=U\otimes_B\kk_\lambda,\qquad
 \nabla(\lambda)=D\Delta(\lambda).
\]
Here $D$ is restricted duality defined using the Chevalley
antiautomorphism; thus $(DM)_\mu=(M_\mu)^*$. The simple head of
$\Delta(\lambda)$ is denoted by $L(\lambda)$.

A \emph{window} is a nonempty finite subset $K$ closed in $J$: if
$\mu\in K$, $\mu\leq\nu$ and $\nu\in J$, then $\nu\in K$.
In particular, $K$ is convex. The subcategory supported in
$J\setminus K$ is Serre, so we may localize by it. We write
\begin{equation}\label{eq:window-category}
 \CC_K=(\OO^J/\OO^{J\setminus K})^{\fl},\qquad
 q_K:\RR_J\longrightarrow\CC_K
\end{equation}
for the finite-length part of the Serre quotient and its exact quotient
functor on $\RR_J$. Its kernel consists of the modules supported in
$J\setminus K$. Put
$\Delta_K(\lambda)=q_K\Delta(\lambda)$,
$\nabla_K(\lambda)=q_K\nabla(\lambda)$ and
$L_K(\lambda)=q_KL(\lambda)$ for $\lambda\in K$.
For $K\subseteq L$, let $q_K^L:\CC_L\to\CC_K$ be the quotient
functor whose kernel has simple labels in $L\setminus K$.
Then $q_K^Lq_L\cong q_K$ and, for $K\subseteq L\subseteq H$,
$q_K^Lq_L^H\cong q_K^H$ canonically.
When all labels have been removed, we use the conventions
$\CC_\varnothing=0$ and $A_\varnothing=0$. The corresponding
dual category and Ringel algebra are also zero.

\subsection{Weight projectives and their endomorphism algebra}
Weight evaluation provides the projectives with which to describe
this localization. For $\mu\in J$ define
\[
 C^J(\mu)=(B\otimes_{U^0}\kk_\mu)^J,\qquad
 Q^J(\mu)=U\otimes_B C^J(\mu).
\]
The superscript on the positive induced module means quotienting by
the sum of its weight spaces outside $J$. Downward closure of $J$
makes this sum a $B$-submodule. The resulting quotient is finite
dimensional, since its weights belong to the finite set
$J_{\geq\mu}$. Fiebig's construction then gives a cyclic projective
$Q^J(\mu)$ with a finite Verma flag and a natural isomorphism
\begin{equation}\label{eq:weight-representer}
 \Hom_{\OO^J}(Q^J(\mu),M)\cong M_\mu
 \qquad(M\in\OO^J);
\end{equation}
see \cite[Theorem~17 and its proof]{Fiebig}.

We assemble all maps between these projectives into the locally
unital algebra
\begin{equation}\label{eq:raw-global}
 A_J=\bigoplus_{\lambda,\mu\in J}
 \Hom_{\OO^J}(Q^J(\mu),Q^J(\lambda)),\qquad
 1_\mu A_J1_\lambda=
 \Hom_{\OO^J}(Q^J(\mu),Q^J(\lambda)).
\end{equation}
Multiplication is opposite composition, so a block
$1_\mu A_J1_\lambda$ acts from the $\lambda$-component to the
$\mu$-component of a left module. Each element of $A_J$ has finite
block support. We use unital modules in the locally unital sense:
a module $V$ satisfies $V=\bigoplus_{\mu\in J}1_\mu V$.

For a window $K$, retain the corresponding idempotents and projectives
by putting
\[
 e_K=\sum_{\mu\in K}1_\mu,\quad A_K=e_KA_Je_K,
 \quad Q_K(\mu)=q_KQ^J(\mu),\quad
 P_K=\bigoplus_{\mu\in K}Q_K(\mu).
\]

\begin{theorem}\label{thm:raw}
The following are equivalences, compatible with the displayed
quotient functor:
\[
 \begin{tikzcd}[column sep=large]
 \RR_J \arrow[r,"\mathcal H_J","\sim"'] \arrow[d,"q_K"'] &
 A_J\text{-}\mathrm{mod}_{\mathrm{lfd}}
       \arrow[d,"e_K(-)"]\\
 \CC_K \arrow[r,"\mathcal H_K","\sim"'] &
 A_K\text{-}\mathrm{mod}_{\mathrm{fd}}.
 \end{tikzcd}
\]
Here $\mathcal H_J(M)=\bigoplus_\mu M_\mu$, and locally finite
dimensional means that every $1_\mu$-component is finite dimensional.
Moreover, applying $q_K$ to the projective maps induces an algebra
isomorphism
\begin{equation}\label{eq:raw-finite}
 A_K\xrightarrow{\sim}
 \End_{\CC_K}(P_K)^{\op},
 \qquad
 \mathcal H_K=\Hom_{\CC_K}(P_K,-).
\end{equation}
The square commutes naturally on objects and morphisms.
\end{theorem}

\begin{proof}
Weight evaluation is exact and commutes with direct sums. Thus the
$Q^J(\mu)$ are small projective generators of $\OO^J$: homogeneous
generators present every object, and then its kernel, by sums of
these projectives. This gives the equivalence with unital $A_J$-modules.
Pointwise finiteness corresponds to finite-dimensional idempotent
components.

Fix $\mu\in K$. In the Serre-quotient Hom formula, the source
subobject $Q'\subseteq Q^J(\mu)$ satisfies
$Q^J(\mu)/Q'\in\OO^{J\setminus K}$. It therefore contains the
distinguished $\mu$-weight generator, so $Q'=Q^J(\mu)$.
The target subobject $N'\subseteq N$ lies in $\OO^{J\setminus K}$,
and projectivity gives
\[
 \Hom(Q^J(\mu),N')=(N')_\mu=0,
 \qquad \Ext^1(Q^J(\mu),N')=0.
\]
Hence replacing $N$ by $N/N'$ also preserves Hom, giving
\begin{equation}\label{eq:retained-projective-hom}
 \Hom_{\OO^J}(Q^J(\mu),N)
 \xrightarrow{\sim}
 \Hom(q_KQ^J(\mu),q_KN).
\end{equation}
This is the retained-weight form of
\cite[Section~17 and Proposition~24]{Fiebig}.

An object has all retained weight spaces zero precisely when it is
supported in $J\setminus K$ and hence is zero in the quotient.
The retained weight functors are therefore exact and jointly faithful.
Each $Q_K(\mu)$ has finite retained dimension, which bounds
its length, and these projectives generate. Every simple is their
quotient, so finite length also implies finite retained dimension.
Thus $P_K$ is a finite projective generator.

Applying \eqref{eq:retained-projective-hom} to the projectives gives
\eqref{eq:raw-finite}. The action is
\[
 a^{\op}\cdot h=h\circ a
\]
for $a\in\End(P_K)$ and $h:P_K\to M$; a map $u:M\to N$
sends $h$ to $u\circ h$. These formulas preserve composition and
identify $q_K$ naturally with $e_K(-)$.
\end{proof}

For $K\subseteq L$, the same realization identifies $q_K^L$ with
$e_K(-)$ and $A_K$ with $e_KA_Le_K$. For every $V\in\CC_L$,
\begin{equation}\label{eq:quotient-end-map}
 \End_{\CC_L}(V)\longrightarrow\End_{\CC_K}(q_K^LV),
 \qquad f\longmapsto q_K^L(f),
\end{equation}
is a unital algebra homomorphism. In the module realization, it
restricts an endomorphism of $\mathcal H_L(V)$ to
$e_K\mathcal H_L(V)$.

\subsection{Global structure and reconstruction}
The weight-projective model gives the supported category an upper
finite highest-weight structure.

\begin{theorem}\label{thm:global-category}
The category $\RR_J$ is upper finite highest weight with poset $J$,
standards $\Delta(\lambda)$ and costandards $\nabla(\lambda)$.
Each stratum is equivalent to $\kk$-\(\mathrm{Vect}_{\mathrm{fd}}\).
Its quotients $\CC_K$ are finite split highest-weight
categories.
\end{theorem}

\begin{proof}
By \eqref{eq:weight-representer}, the blocks of $A_J$ are weight
spaces of the pointwise finite projectives $Q^J(\lambda)$, and are
therefore finite dimensional. Theorem~\ref{thm:raw} realizes
$\RR_J$ as a Schurian category in the sense of
\cite[Section~2.3]{BS}, with upper finite set of simple labels $J$.
We first check that downward support conditions agree with the
corresponding conditions on simple subquotient labels. A nonzero
vector of weight $\mu$ generates a nonzero cyclic quotient of
$Q^J(\mu)$. Its head is nonzero and finite by \cite[Lemma~2.14]{BS}.
The generator survives in some simple quotient $L(\nu)$, so
$\nu\geq\mu$. Downward closure then recovers the weight condition
from the condition on simple labels.

Inside the category supported below $\lambda$, evaluation at
$\lambda$ is represented by $\Delta(\lambda)$, with endomorphism
field $\kk$. Write $\RR_{\leq\lambda}$ and $\RR_{<\lambda}$ for
the full subcategories supported on the indicated sets of weights.
The stratum $\RR_{\leq\lambda}/\RR_{<\lambda}$ is therefore
$\kk$-\(\mathrm{Vect}_{\mathrm{fd}}\), and restricted duality
supplies the costandard object. To verify the projective filtration
axiom, use the positive inducing module to obtain
\begin{equation}\label{eq:standard-projective-sequence}
 0\longrightarrow F^J(\lambda)\longrightarrow Q^J(\lambda)
 \longrightarrow\Delta(\lambda)\longrightarrow0,
\end{equation}
where $F^J(\lambda)$ has a finite Verma flag with labels strictly
larger than $\lambda$. These are the projective filtration axioms
of \cite[Theorem~3.33 and Definition~3.34]{BS}. Applying the quotient
functors $q_K$ gives the finite highest-weight categories.
\end{proof}

Because the strata are simple, proper standards coincide with
standards and are unaffected by the choice of signs in a
stratification. Thus $\RR_J$ is also fully stratified.

Compatible finite images recover both objects and morphisms.

\begin{proposition}\label{prop:coherent-reconstruction}
The functor $M\mapsto(q_KM)_K$ identifies $\RR_J$ with the category
of coherent finite-window systems. Such a system consists of
$M_K\in\CC_K$ and isomorphisms
$\theta_{L,K}:q_K^LM_L\to M_K$. They satisfy $\theta_{K,K}=1$ and,
for $K\subseteq L\subseteq H$,
\[
 \theta_{H,K}=\theta_{L,K}\circ q_K^L(\theta_{H,L}),
\]
using the canonical identification $q_K^Lq_L^H\cong q_K^H$.
Thus passing directly from $H$ to $K$ agrees with passing through $L$.
Morphisms are families $f_K:M_K\to N_K$ commuting with these
identifications. In particular,
\begin{equation}\label{eq:global-hom-limit}
 \Hom_U(M,N)\xrightarrow{\sim}
 \varprojlim_K\Hom_{\CC_K}(q_KM,q_KN)
 \qquad(M,N\in\RR_J),
\end{equation}
where the transition from $L$ to $K$ sends $f$ to $q_K^L(f)$.
These isomorphisms preserve composition.
\end{proposition}

\begin{proof}
In the weight-projective model a coherent system is a family of finite
$A_K$-modules $V_K$ with $e_KV_L\cong V_K$. For each $\mu\in J$,
the given identifications identify the $\mu$-components of all
windows containing $\mu$ with one finite-dimensional space $V_\mu$.
Set $V=\bigoplus_\mu V_\mu$. To define the action of a block
$a\in1_\alpha A_J1_\beta$, choose a finite closed subset containing
$\alpha$ and $\beta$ and use its action on the corresponding
components. Such a closed subset exists by quasi-boundedness, and
coherence makes the result independent of the choice. For composable
blocks in $1_\alpha A_J1_\beta$ and $1_\beta A_J1_\gamma$,
choose one window containing $\alpha,\beta,\gamma$; associativity
there verifies the product on $V$. Consequently $V$ is a
unital, locally finite-dimensional $A_J$-module and corresponds to
an object of $\RR_J$. Compatible local maps assemble componentwise
by the same argument; their compositions are also computed in a
common window. This proves the equivalence and
\eqref{eq:global-hom-limit}, including compatibility with composition.
\end{proof}

The $A_J$-module realization gives the left and right adjoints
$j_{K,!},j_{K,*}:\CC_K\to\RR_J$ of $q_K$ as induction and
coinduction, respectively \cite[Section~2.5]{BS}:
\begin{equation}\label{eq:window-adjoints}
 j_{K,!}V=A_Je_K\otimes_{A_K}V,
 \qquad
 (j_{K,*}V)_\mu=\Hom_{A_K}(e_KA_J1_\mu,V).
\end{equation}
Both functors are fully faithful. Their actions on morphisms
are tensoring and postcomposition, respectively. Induction is right
exact and coinduction is left exact. We will use their stronger
exactness properties on the appropriate flag classes, as provided
by \cite[Theorem~3.42]{BS}.

\subsection{Frobenius realization}\label{sec:frobenius}

Frobenius restriction retains the PBW operators whose divided-power
exponents are less than $p^r$. The window width determines when these
already give the whole algebra.

\subsubsection{Restriction and bounded PBW exponents}
Let $r\geq1$, set $q=p^r$ and write $H_r=\Dist(G_r)$. We retain
the rational torus and define
\[
 \CC_{r,K}=
 \bigl(\operatorname{Rep}(G_rT)^J/
       \operatorname{Rep}(G_rT)^{J\setminus K}\bigr)^{\fl}.
\]
Here $\operatorname{Rep}(G_rT)$ denotes all rational $G_rT$-modules.
Write $q_{r,K}$ for the quotient functor, restricted to modules with
finite-dimensional weight spaces, and
$q_{r,K}^L:\CC_{r,L}\to\CC_{r,K}$ for its counterpart when
$K\subseteq L$. Their kernels are given by discarded support.
Restriction preserves the two support subcategories, so
it descends through their Serre quotients. The weight-projective
construction below shows that the image of a finite-length window
object again has finite length, giving the functor
$R_{r,K}:\CC_K\to\CC_{r,K}$.

Put $H_r^{\geq0}=H_r^0H_r^+$, and let $C_r^J(\lambda)$ be the quotient of
$H_r^+\otimes\kk_\lambda$ by its weights outside $J$, and set
\[
 Q_r^J(\lambda)=H_r\otimes_{H_r^{\geq0}}C_r^J(\lambda).
\]
The retained torus acts by
$t(h\otimes c)=(\operatorname{Ad}t)(h)\otimes tc$.
Adjunction gives
$\Hom_{G_rT}(Q_r^J(\lambda),N)=N_\lambda$.
The localization argument in
\eqref{eq:retained-projective-hom} applies to these cyclic
projectives as well. Their quotient images, summed over
$\lambda\in K$, therefore form a projective generator
$P_{r,K}$. Put
$A_{r,K}=\End(P_{r,K})^{\op}$.

Fix orders on $R^+$ and write
$f^{(a)}=\prod_{\beta\in R^+}f_\beta^{(a_\beta)}$ and
$e^{(c)}=\prod_{\beta\in R^+}e_\beta^{(c_\beta)}$ in these orders.
For an exponent vector $a$, put
$|a|=\sum_{\beta\in R^+}a_\beta\beta$.
A basis of $1_\mu A_K1_\lambda$ is
\begin{equation}\label{eq:peak-basis}
 f^{(a)}e^{(c)}1_\lambda,
 \qquad \nu=\lambda+|c|=\mu+|a|\in K.
\end{equation}
Here $\nu$ is the intermediate weight, or peak, of the monomial.
The corresponding basis for $1_\mu A_{r,K}1_\lambda$ has the
same form, subject to the additional restrictions
\begin{equation}\label{eq:capped-pbw}
 0\leq a_\beta,c_\beta<q\qquad(\beta\in R^+).
\end{equation}
These descriptions follow from the PBW bases of the two induced
weight projectives. In the positive factor the intermediate weight
must lie in $J$; closedness of $K$ in $J$ then places it in $K$, since
it lies above a retained weight. Thus the two bases use precisely
the same allowed intermediate weights. Their only difference is
the bound on each divided-power exponent.

\begin{proposition}\label{prop:actual-inclusion}
The PBW operators satisfying \eqref{eq:capped-pbw} define a unital
inclusion $\iota_{r,K}:A_{r,K}\hookrightarrow A_K$. Under the projective
generator equivalences, $R_{r,K}$ is restriction of scalars along
$\iota_{r,K}$, with action $a\cdot v=\iota_{r,K}(a)v$.
In particular, $R_{r,K}$ is exact and faithful for every $r\geq1$.
\end{proposition}

\begin{proof}
The two projective constructions establish independence of their
respective PBW bases, so sending a basis element satisfying
\eqref{eq:capped-pbw} to the
same operator in $A_K$ is injective. This map also preserves
multiplication: $H_r\subseteq U$ is a subalgebra, toral distributions
act by the fixed integral weight scalars, and both positive
inductions impose the same outside-$J$ relations. It preserves
the identity $e_K=\sum_{\lambda\in K}1_\lambda$.

On a restricted module, each of these operators acts by the same
distribution as before restriction. Its retained-weight
representation is therefore restriction of scalars along the
displayed inclusion. This description proves exactness and
faithfulness, and also the finite-length assertion, since the
retained representation remains finite dimensional.
\end{proof}

\subsubsection{The exact bound}
For a comparable pair $\lambda\leq\mu$ in $K$, write
$\mu-\lambda=\sum_i n_i\alpha_i$. Define the width
\[
 d(K)=\max\bigl(\{n_i:\lambda,\mu\in K,
                          \ \lambda\leq\mu\}\cup\{0\}\bigr).
\]

\begin{theorem}\label{thm:frobenius}
The restriction functor $R_{r,K}$ is an equivalence if and
only if $p^r>d(K)$.
\end{theorem}

\begin{proof}
Suppose first that $q>d(K)$. Any PBW exponent at least $q$ in
\eqref{eq:peak-basis} would force some simple-root coefficient
of the difference from an endpoint to $\nu$ to be at least $q$.
Both endpoints and $\nu$ belong to $K$, so this contradicts the
definition of $d(K)$. Every basis element therefore satisfies
\eqref{eq:capped-pbw}: the inclusion in Proposition~\ref{prop:actual-inclusion}
is an isomorphism, and its restriction-of-scalars functor is an
equivalence.

Conversely, suppose $q\leq d(K)$. Choose comparable endpoints
whose difference has an $\alpha_i$-coefficient at least $q$.
Taking the lower endpoint to be $\lambda$, the weight
$\nu=\lambda+q\alpha_i$ lies between the two endpoints in the root
order. Convexity therefore gives $\nu\in K$.
Let $\mathcal P(\eta)$ be the positive-root partition function,
and $\mathcal P_q(\eta)$ its version with each exponent less
than $q$. The two representing projectives at $\lambda$ have
endomorphism dimensions
\[
 \sum_{\xi\in K}\mathcal P(\xi-\lambda)^2,
 \qquad
 \sum_{\xi\in K}\mathcal P_q(\xi-\lambda)^2.
\]
Every summand in the first sum is at least the corresponding
summand in the second. At $\xi=\nu$ the inequality is strict:
the only positive-root partition of $q\alpha_i$ consists of
$q$ copies of $\alpha_i$, so
$\mathcal P(q\alpha_i)=1$ and $\mathcal P_q(q\alpha_i)=0$.
If $R_{r,K}$ were an equivalence, its preservation of the
$\lambda$-weight functor would identify the two representing
projectives by Yoneda's lemma. Their endomorphism algebras would then
have the same dimension, contradicting the strict inequality.
\end{proof}

When the bound holds, the equivalence transports the highest-weight
structure to $\CC_{r,K}$. A selected family may be recovered earlier.

These equivalences respect passage between windows. For
$K\subseteq L$, the support quotients give the natural square
\begin{equation}\label{eq:nested-frobenius-square}
 \begin{tikzcd}[column sep=large]
 \CC_L \arrow[r,"q_K^L"] \arrow[d,"R_{r,L}"'] &
 \CC_K \arrow[d,"R_{r,K}"]\\
 \CC_{r,L} \arrow[r,"q_{r,K}^L"'] & \CC_{r,K}.
 \end{tikzcd}
\end{equation}
Its commutativity follows by restricting representatives and their
localization maps; both vertical arrows are equivalences when
$p^r>d(L)$. The same construction compares Frobenius levels:
for $s\geq r$, restriction from $G_sT$ to $G_rT$ induces
$H_{s,r,K}:\CC_{s,K}\to\CC_{r,K}$, and
$H_{s,r,K}R_{s,K}\cong R_{r,K}$. If $s\geq r$ and
$p^r>d(K)$, this is an exact equivalence. In particular it
preserves all quotient Yoneda Ext groups and their products.

\subsubsection{Periodicity and the retained torus}
The retained torus also records the weight shift in Fiebig's
periodicity.

For this paragraph use the canonical support
$J=\operatorname{down}(K)$. Choose a dominant weight $\gamma$
and an integer $\ell\geq1$ with $p^\ell>d(K)$, and put
$\tau=p^\ell\gamma$. The width bound ensures that $K$ contains
no pair differing by $p^\ell\eta$ with $0\ne\eta\in\NN\Pi$.
Tensoring by
$L(\tau)=L(\gamma)^{[\ell]}$ induces Fiebig's equivalence
$E_{\gamma,\ell}:\CC_K\to\CC_{K+\tau}$
\cite[Theorem~25]{Fiebig}. If $1\leq r\leq\ell$, there is a
natural commutative square
\begin{equation}\label{eq:periodic-square}
 \begin{tikzcd}[column sep=large]
 \CC_K \arrow[r,"E_{\gamma,\ell}","\sim"']
       \arrow[d,"R_{r,K}"'] &
 \CC_{K+\tau} \arrow[d,"R_{r,K+\tau}"]\\
 \CC_{r,K} \arrow[r,"-\otimes\kk_\tau","\sim"'] &
 \CC_{r,K+\tau}.
 \end{tikzcd}
\end{equation}
To verify the square, choose a nonzero highest-weight vector in
$L(\tau)$. Since $r\leq\ell$, the $G_r$-action on this Frobenius
twist is trivial, while its highest torus line retains the weight
$\tau\in p^rX$. The inclusion $\kk_\tau\hookrightarrow L(\tau)$
therefore gives a natural map
$\operatorname{Res}_{G_rT}M\otimes\kk_\tau
\to\operatorname{Res}_{G_rT}(M\otimes L(\tau))$.
We claim that every lower line disappears in
the shifted quotient. Such a line has weight
$\tau-p^\ell\eta$ with
$0\ne\eta\in\NN\Pi$. If its tensor product with a supported
module had a retained weight $\kappa+\tau$, then some
$\mu\in J$ would satisfy $\kappa=\mu-p^\ell\eta$. Choose
$\xi\in K$ above $\mu$. The difference $\xi-\kappa$ has a
simple-root coefficient at least $p^\ell$, contradicting the
width bound. Thus only the highest line survives the quotient,
and tensoring with it gives the lower functor in the square.
The resulting map respects discarded support and is natural
before and after localization, proving the claimed compatibility.

The vertical functors are equivalences when $p^r>d(K)$.
% END foundations.tex

% BEGIN tiltings.tex
\section{Global tiltings and finite recovery}\label{sec:tiltings}

We retain the category $\RR_J=\OO^J_{\pf}$ and the finite closed subsets of
the fixed support $J$ from Section~\ref{sec:windows}. A finite standard
flag is a finite filtration whose sections are ordinary Verma modules;
a finite costandard flag uses their restricted duals. To allow infinite
flags, we call an exhaustive directed union of finite standard-filtered
subobjects $M_a$ an \emph{ascending standard flag} if $M_b/M_a$ has a
finite standard flag whenever $a\leq b$. Descending costandard flags
are obtained by restricted duality. A \emph{global tilting} is an object
of $\RR_J$ with both an ascending standard flag and a descending
costandard flag.

The finite quotients detect tiltings; the two adjoints determine
when a single quotient recovers all their outgoing morphisms.

\subsection{Tiltings and their finite images}

In the finite highest-weight category $\CC_K$, let $T_K(\lambda)$ be the
indecomposable tilting with highest label $\lambda\in K$. Its existence
and uniqueness are the finite highest-weight case of
\cite[Theorem~4.2]{BS}. The same theorem gives
$(T_K(\lambda):\Delta_K(\lambda))=1$, with every other standard label
strictly below $\lambda$. Restricted duality fixes $T_K(\lambda)$.

For finite standard-filtered $V$ and costandard-filtered $W$ in
$\CC_K$, standard--costandard orthogonality gives
\begin{equation}\label{eq:finite-flag-hom}
 \dim\Hom_{\CC_K}(V,W)
   =\sum_{\nu\in K}(V:\Delta_K(\nu))(W:\nabla_K(\nu)),
 \qquad \Ext^1_{\CC_K}(V,W)=0.
\end{equation}
Indeed, the corresponding assertions for individual standard and
costandard objects are
$\Hom(\Delta_K(\nu),\nabla_K(\eta))=\delta_{\nu\eta}\kk$ and
$\Ext^1(\Delta_K(\nu),\nabla_K(\eta))=0$; induction on the two finite
flags proves the formula.

To pass between finite and global flags, we use the left and right
adjoints $j_{K,!},j_{K,*}$ of $q_K$ in
\eqref{eq:window-adjoints}. The left adjoint is exact on standard flags,
and the right adjoint is exact on costandard flags; see
\cite[Theorem~3.42(3)--(6)]{BS}. We compute $\Ext_{\OO^J}$ in the
category of all unital $A_J$-modules, using the realization in
Theorem~\ref{thm:raw} and the convention of \cite[Section~2.3]{BS}.

\begin{lemma}\label{lem:flag-detection}
An object $M\in\RR_J$ has an ascending standard flag, respectively a
descending costandard flag, if and only if every $q_KM$ has a finite
standard flag, respectively a finite costandard flag. Thus $M$ is a
global tilting if and only if all its finite images are tilting.
\end{lemma}
\begin{proof}
For $\nu\in K$, the derived adjunctions of
\cite[Theorem~3.42(3)--(6)]{BS} give, for every $n\geq0$,
\begin{align}
 \Ext^n_{\OO^J}(\Delta(\nu),M)
   &\cong\Ext^n_{\CC_K}(\Delta_K(\nu),q_KM),\label{eq:standard-adjunction}\\
 \Ext^n_{\OO^J}(M,\nabla(\nu))
   &\cong\Ext^n_{\CC_K}(q_KM,\nabla_K(\nu)).\label{eq:costandard-adjunction}
\end{align}
If every finite image has a standard flag, the second isomorphism
gives vanishing against all costandards, which detects an ascending
standard flag by \cite[Theorem~3.37]{BS}. If the finite images have
costandard flags, the first gives vanishing from all standards, which
detects a descending costandard flag by \cite[Theorem~3.39]{BS}.
For the other direction, applying the
exact functor $e_K(-)$ to an ascending flag gives a directed union of
standard-filtered submodules of the finite-dimensional module $e_KM$.
This union stabilizes at one member, so the local flag is finite.
Restricted duality gives the costandard assertion.
\end{proof}

For $K\subseteq L$, the quotient functor
$q_K^L:\CC_L\to\CC_K$ induces a surjection
\begin{equation}\label{eq:finite-hom-lift}
 \Hom_{\CC_L}(V,W)\twoheadrightarrow
 \Hom_{\CC_K}(q_K^LV,q_K^LW)
\end{equation}
whenever $V$ is standard filtered and $W$ is costandard filtered;
the map sends $f$ to $q_K^L(f)$.
To see this, work in the weight-projective model and set
$V^-=A_Le_KV$, the submodule generated by the $K$-weight spaces.
Its standard-flag labels lie in $K$. Let $W^+$ be the largest subobject of $W$ supported
in $L\setminus K$. The finite truncation theorem
\cite[Theorem~3.18 and Corollary~3.19]{BS} identifies
$V^-\cong j_!q_K^LV$ and $W/W^+\cong j_*q_K^LW$.
Here $j_!$ and $j_*$ are its left and right adjoints, respectively;
both $V^-$ and
$V/V^-$ are standard filtered, while $W^+$ is costandard filtered.
Adjunction identifies a local map with a map $V^-\to W/W^+$.
Applying $\Hom(V^-,-)$ to $0\to W^+\to W\to W/W^+\to0$ lifts
this map to $V^-\to W$, since
$\Ext^1_{\CC_L}(V^-,W^+)=0$.
Applying $\Hom(-,W)$ to $0\to V^-\to V\to V/V^-\to0$ then
extends it to $V\to W$, since
$\Ext^1_{\CC_L}(V/V^-,W)=0$.
Both vanishings follow from \eqref{eq:finite-flag-hom}, applied in
$\CC_L$. The resulting map has the prescribed image under $q_K^L$,
which proves surjectivity. Functoriality preserves composition.

\begin{theorem}\label{thm:tilting}
For every $\lambda\in J$ there is a unique, up to isomorphism,
indecomposable global tilting $T(\lambda)\in\RR_J$
with highest label $\lambda$ and highest weight space of dimension one.
It is restricted self-dual and satisfies
\[
 q_KT(\lambda)\cong
 \begin{cases}T_K(\lambda),&\lambda\in K,\\0,&\lambda\notin K.
 \end{cases}
\]
For every $\nu\in J$, the groups
$\Ext^n_{\OO^J}(\Delta(\nu),T(\lambda))$ and
$\Ext^n_{\OO^J}(T(\lambda),\nabla(\nu))$ vanish for
$n>0$.
\end{theorem}
\begin{proof}
We assemble the finite tiltings and identify the result by the
uniqueness theorem of \cite[Theorem~4.18]{BS}. Choose finite closed subsets
$K_0\subseteq K_1\subseteq\cdots$ containing $\lambda$ with union $J$.
The sequence is \emph{cofinal}, meaning that every window lies in
some $K_n$. Take upper closures in $J$ of finite initial segments
of an enumeration of $J$; upper finiteness makes each closure finite.
Tilting truncation \cite[Corollary~4.5]{BS}
identifies each quotient of $T_{K_{n+1}}(\lambda)$ with
$T_{K_n}(\lambda)$. Composing successive identifications makes them
coherent, and Proposition~\ref{prop:coherent-reconstruction} glues them
to a pointwise finite object. Truncation gives the stated images in
every window, and Lemma~\ref{lem:flag-detection} gives both flags.

An idempotent endomorphism is zero or the identity in each finite
indecomposable tilting. Coherent unital transitions force the same
choice throughout, so Hom reconstruction proves indecomposability.
The top weight has dimension one. Uniqueness identifies the result
with the asserted tilting and also with its restricted dual.
Finally, choose a window retaining $\lambda,\nu$ in
\eqref{eq:standard-adjunction}--\eqref{eq:costandard-adjunction};
finite standard--costandard orthogonality in every positive degree
gives the Ext vanishings.
\end{proof}

If $K$ retains $\lambda$ and $\nu$, put
\begin{equation}\label{eq:stable-multiplicity}
 t_\lambda(\nu)=(T_K(\lambda):\Delta_K(\nu)),
 \qquad S(T(\lambda))=\{\nu\in J:t_\lambda(\nu)>0\}.
\end{equation}
These integers are independent of $K$: nested quotients preserve
standards at retained labels and kill the others, and two windows
can be compared inside a common larger window. They vanish unless
$\nu\leq\lambda$. We call $S(T)$ the \emph{standard support}.

\begin{proposition}\label{prop:tilting-hom}
Let $\lambda,\mu\in J$ and put
$S_{\lambda,\mu}=\{\nu:t_\lambda(\nu)t_\mu(\nu)>0\}$.
The maps induced by $q_K$ give
\begin{equation}\label{eq:tilting-hom-limit}
 \Hom_U(T(\lambda),T(\mu))
 \cong\varprojlim_K\Hom_{\CC_K}(T_K(\lambda),T_K(\mu)),
\end{equation}
where the windows retain both labels. Write $\rho_K(f)=q_K(f)$
for the induced map to the $K$-term. It is surjective,
and is an isomorphism precisely when $S_{\lambda,\mu}\subseteq K$.
The global Hom space is finite dimensional precisely when
$S_{\lambda,\mu}$ is finite; in that case
\[
 \dim\Hom_U(T(\lambda),T(\mu))
   =\sum_{\nu\in J}t_\lambda(\nu)t_\mu(\nu),\qquad
 \dim\ker\rho_K
   =\sum_{\nu\notin K}t_\lambda(\nu)t_\mu(\nu).
\]
If $S_{\lambda,\mu}$ is infinite, every $\rho_K$ has an
infinite-dimensional kernel.
\end{proposition}
\begin{proof}
The isomorphism \eqref{eq:tilting-hom-limit} follows from
Proposition~\ref{prop:coherent-reconstruction}: a compatible family
of maps is determined by its action on every weight component.
Under \eqref{eq:tilting-hom-limit}, $\rho_K$ is the projection to
the $K$-coordinate. To lift a specified map at $K$, start a cofinal
sequence at $K$ and lift successively using
\eqref{eq:finite-hom-lift}. The resulting compatible family defines
a global map by Proposition~\ref{prop:coherent-reconstruction}.

The local dimensions are
$\sum_{\nu\in K}t_\lambda(\nu)t_\mu(\nu)$ by
\eqref{eq:finite-flag-hom} and self-duality. When the common support
is finite, these dimensions stabilize, and surjectivity makes all
later transitions isomorphisms. The inverse limit is then identified
with any sufficiently large term, proving both dimension formulas.
If the support is infinite, the finite dimensions are unbounded, so the
global space is infinite dimensional. At a fixed window its target
is finite dimensional, and hence its kernel is infinite dimensional.
Finally, if $\nu\in S_{\lambda,\mu}\setminus K$, choose $L\supseteq K$
containing $\nu$. The local dimension strictly increases, so the
surjective transition from $L$ to $K$ has a nonzero kernel.
Surjectivity of $\rho_L$ lifts a nonzero element of this kernel
globally, proving that $\rho_K$ is not injective.
\end{proof}

\subsection{Canonical lifts and the finiteness criterion}

To recover maps to arbitrary targets, we use the canonical lifts of a
finite image.

\begin{lemma}\label{lem:canonical-lifts}
If $M\in\RR_J$ has an ascending standard flag, the adjunction counit
identifies $j_{K,!}q_KM$ with
\[
 M_K^-=U\sum_{\nu\in K}M_\nu\ \subseteq M.
\]
This subobject has a finite standard flag with the retained
multiplicities, and $M=\bigcup_KM_K^-$. If $M$ has a descending
costandard flag, the unit identifies $j_{K,*}q_KM$ with $M/M_K^+$,
where $M_K^+$ is the largest submodule supported in $J\setminus K$.
The quotient has a finite costandard flag with the retained
multiplicities.
\end{lemma}
\begin{proof}
In the module model the counit is the multiplication map
$A_Je_K\otimes_{A_K}e_KM\to M$, whose image is $M_K^-$.
In every finite window containing $K$, \cite[Corollary~3.19]{BS}
identifies the counit with the inclusion of the submodule generated
by the retained weight spaces.
The kernel therefore has zero image in every such window, hence
every weight component is zero by
Proposition~\ref{prop:coherent-reconstruction}. Thus the counit is
monic. Its finite flag
follows by applying $j_{K,!}$ to the local flag, using
\cite[Theorem~3.42]{BS}. Every weight vector belongs to some retained
window, proving the union. Restricted duality gives the corresponding
epimorphism $M\to j_{K,*}q_KM$. Its kernel is killed by $q_K$;
conversely, adjunction forces every submodule killed by $q_K$ to map
to zero under the unit. This proves the asserted maximality of
$M_K^+$ and completes the costandard statement.
\end{proof}

\begin{theorem}\label{thm:recovery}
For $T=T(\lambda)$ the following conditions are equivalent:
\begin{enumerate}
\item $T$ is finitely generated.
\item $T$ is finitely presented in $\OO^J$.
\item $T$ has finite ordinary standard and costandard flags.
\item $S(T)$ is finite.
\item $\End_U(T)$ is finite dimensional.
\item For some window $K$ retaining $\lambda$, the map induced by $q_K$,
      $\End_U(T)\to\End_{\CC_K}(T_K(\lambda))$ is an isomorphism.
\item For some such $K$, the induced map
\[
 \Hom_U(T,N)\longrightarrow\Hom_{\CC_K}(q_KT,q_KN)
\]
is an isomorphism for every $N\in\RR_J$.
\end{enumerate}
In the last condition one may take any finite closed subset containing
$S(T)$.
\end{theorem}
\begin{proof}
We first relate finite generation to the stable multiplicities.
Evaluation at finitely many homogeneous generators embeds $\End_U(T)$
in a finite sum of finite-dimensional weight spaces. Thus (i) implies
(v). Proposition~\ref{prop:tilting-hom}, with $\mu=\lambda$, proves
(iv)$\Leftrightarrow$(v)$\Leftrightarrow$(vi).

Now suppose $S(T)\subseteq K$. For every larger finite window $L$, the
submodule of $q_LT$ generated by the $K$-weight spaces has the same
standard multiplicities as $q_LT$. The quotient by this submodule
is standard filtered with no standard sections, and is therefore
zero. Proposition~\ref{prop:coherent-reconstruction} and
Lemma~\ref{lem:canonical-lifts} give
$T=T_K^-=j_{K,!}q_KT$. Thus $T$ has a finite standard flag; self-duality
gives a finite costandard flag. The converse (iii)$\Rightarrow$(iv)
follows by applying the quotients to a finite flag.

To complete the finiteness implications, each $\Delta(\nu)$ is
finitely presented. Indeed, the projective
sequence
$0\to F^J(\nu)\to Q^J(\nu)\to\Delta(\nu)\to0$
from Fiebig's construction has a finite standard flag on its kernel
\cite[proof of Theorem~17]{Fiebig}. Standards are cyclic, so this kernel is
finitely generated and is a quotient of a finite sum of weight
projectives. This gives a finite presentation. Finite extensions of
finitely presented modules are finitely presented, proving
(iii)$\Rightarrow$(ii); (ii)$\Rightarrow$(i) is immediate.

Finally, whenever $S(T)\subseteq K$, the canonical identification
$T\cong j_{K,!}q_KT$ makes the map in (vii) the adjunction
isomorphism, for every $N$. Taking $N=T$ proves
(vii)$\Rightarrow$(vi).
\end{proof}

For $V\in\CC_K$, let
$c_{K,V}:j_{K,!}V\to j_{K,*}V$ be the map corresponding to
$\mathrm{id}_V$ under adjunction. It is natural in $V$, and $q_Kc_{K,V}$
is the identity. For $T=T(\lambda)$ and $\lambda\in K$, write
$c_K=c_{K,T_K(\lambda)}$. Lemma~\ref{lem:canonical-lifts} identifies
it with the composite $T_K^-\hookrightarrow T\twoheadrightarrow
T/T_K^+$. Consequently
\begin{equation}\label{eq:comparison-defects}
\begin{split}
 0\longrightarrow T_K^-\cap T_K^+
 \longrightarrow j_{K,!}T_K(\lambda)
 \xrightarrow{\ c_K\ }j_{K,*}T_K(\lambda)\\
 \longrightarrow T/(T_K^-+T_K^+)
 \longrightarrow0
\end{split}
\end{equation}
is exact: the kernel is $T_K^-\cap T_K^+$ and the cokernel is
$T/(T_K^-+T_K^+)$, both supported in $J\setminus K$.

\begin{theorem}\label{thm:adjoint-comparison}
For a finite closed subset $K$ retaining $\lambda$, the following are
equivalent: $S(T)\subseteq K$; $c_K$ is an isomorphism; $c_K$ is a
monomorphism; $c_K$ is an epimorphism. When they hold, both canonical
maps
\[
 j_{K,!}T_K(\lambda)\longrightarrow T
       \longrightarrow j_{K,*}T_K(\lambda)
\]
are isomorphisms. If $S(T)$ is finite, the smallest such closed subset is
$\uparrow_JS(T)=\{\nu\in J:\mu\leq\nu\text{ for some }\mu\in S(T)\}$.
\end{theorem}
\begin{proof}
If $S(T)\subseteq K$, the preceding proof gives $T=T_K^-$, and
self-duality gives $T_K^+=0$. Conversely, if $c_K$ is an isomorphism,
\eqref{eq:comparison-defects} gives
$T_K^-\cap T_K^+=0$ and $T_K^-+T_K^+=T$, hence
$T=T_K^-\oplus T_K^+$. Since $T$ is indecomposable and $q_KT\ne0$,
the second summand is zero. Thus $T=T_K^-$ has a finite standard flag
with all labels in $K$.

The two lifts have equal characters:
\[
 \operatorname{ch}(j_{K,!}T_K(\lambda))
  =\sum_{\nu\in K}t_\lambda(\nu)\operatorname{ch}\Delta(\nu)
  =\operatorname{ch}(j_{K,*}T_K(\lambda)).
\]
Here self-duality gives equal local standard and costandard
multiplicities. Their weight spaces are finite dimensional, so a
monomorphism or epimorphism between these particular objects is an
isomorphism. Finally, upper finiteness makes $\uparrow_JS(T)$ finite
when $S(T)$ is finite, and every closed subset containing $S(T)$ contains
this upper closure.
\end{proof}

In particular, if $S(T)$ is infinite, every finite comparison $c_K$
has both a nonzero kernel and a nonzero cokernel. At a one-label
window, with $J=\operatorname{down}(\lambda)$, this comparison is the
usual map $\Delta(\lambda)\to\nabla(\lambda)$ with image $L(\lambda)$.
% END tiltings.tex

% BEGIN family-recovery.tex
\subsection{Recovery of a finite tilting family}
The support criterion also controls finite families. Write $\add(M)$
for the category of direct summands of finite direct sums of $M$.

\begin{corollary}\label{cor:family-quotient}
Let $\Gamma\subseteq J$ be finite, put
$T_\Gamma=\bigoplus_{\lambda\in\Gamma}T(\lambda)$, and suppose that
$K\subseteq J$ is a finite closed subset containing $\Gamma$. The following
conditions are equivalent, with each displayed map induced by $q_K$:
\begin{enumerate}
\item $\End_U(T_\Gamma)\to\End_{\CC_K}(q_KT_\Gamma)$ is an isomorphism;
\item $\End_U(T(\lambda))\to\End_{\CC_K}(T_K(\lambda))$ is an
isomorphism for every $\lambda\in\Gamma$;
\item $\bigcup_{\lambda\in\Gamma}S(T(\lambda))\subseteq K$.
\end{enumerate}
They are also equivalent to full faithfulness of $q_K$ on
$\add(T_\Gamma)$.
\end{corollary}
\begin{proof}
The diagonal idempotent corners of the map in (i) are the individual
endomorphism maps in (ii), so (i) implies (ii).
The individual support criterion proves the equivalence of (ii) and
(iii). Under (iii), the common standard support of every ordered pair
$T(\lambda),T(\mu)$ is retained. The mixed Hom formula and its quotient
map therefore show that every block of the map in (i) is an
isomorphism. These maps are induced by $q_K$, so their
direct sum respects multiplication. Passing to matrices and idempotent
summands proves the last assertion.
\end{proof}

% END family-recovery.tex

% BEGIN ringel.tex
\section{Global Ringel duality and derived realization}\label{sec:ringel}

Semi-infinite Ringel duality turns the quotients of $\RR_J$ into
finite Serre subcategories. Their algebras and transition maps are
induced from the same global tiltings.

Put
\[
 \TJ=\bigoplus_{\lambda\in J}T(\lambda),\qquad
 \Theta_K=q_K\TJ=\bigoplus_{\lambda\in K}T_K(\lambda),\qquad
 \BK=\End_{\CC_K}(\Theta_K)^{\op}.
\]
Quasi-boundedness ensures that the sum $\TJ$ belongs to $\RR_J$:
at a fixed weight $\mu$, only the finitely many labels
$\lambda\geq\mu$ can contribute. Write
$\varepsilon_{\lambda,K}\in\BK$ for the endomorphism that is the
identity on $T_K(\lambda)$ and zero on the other summands, regarded
in the opposite algebra. The sum $\TJ$ is the tilting generator of
\cite[Definition~4.26]{BS}.

Set $\BJ=\End_U(\TJ)^{\op}$. The quotient functors define
\[
 \pi_K:\BJ\longrightarrow\BK,\quad f^{\op}\longmapsto(q_Kf)^{\op},
 \qquad
 \pi_{L,K}:\mathsf B_L\longrightarrow\BK,\quad
 f^{\op}\longmapsto(q_K^Lf)^{\op}.
\]
We use the chosen identifications $q_K\TJ=\Theta_K$ and
$q_K^L\Theta_L=\Theta_K$. Give $\kk$ and each $\BK$ the discrete
topology, and give $\BJ$ the topology with
$I_K=\ker\pi_K$ as a neighbourhood basis at zero.
Thus $b+I_K$ consists of the elements with the same image as $b$
in $\BK$. Since $I_L\subseteq I_K$ for $K\subseteq L$ and a
finite union of windows is a window, these two-sided ideals define
a topological algebra. The next theorem identifies this topology
with the profinite topology of \cite[Section~2.1, Lemma~2.10]{BS}.

\begin{theorem}\label{thm:ringel}
The homomorphisms $\pi_{L,K}$ are surjective and unital, with
\begin{equation}\label{eq:ringel-kernel}
 \ker\pi_{L,K}=
 \mathsf B_L\left(\sum_{\nu\in L\setminus K}
           \varepsilon_{\nu,L}\right)\mathsf B_L.
\end{equation}
There is an isomorphism of topological algebras
\begin{equation}\label{eq:ringel-limit}
 \BJ\xrightarrow{\sim}
       \varprojlim_K\BK.
\end{equation}
Let
\[
 \DD_J:=\BJ\text{-}\mathrm{mod}_{\mathrm{fd,cont}}
\]
be the category of finite-dimensional left modules given the discrete
topology and a continuous $\BJ$-action. Equivalently, their action
factors through $\BK$ for some finite closed subset $K$.
This category is lower finite highest weight with label order $J^{\op}$.
Inflation along $\pi_K:\BJ\twoheadrightarrow\BK$ gives an exact fully
faithful functor
$i_K:\DD_K:=\BK\text{-}\mathrm{mod}_{\mathrm{fd}}\to\DD_J$.
Its essential image is the Serre subcategory
of $\DD_J$ whose simple composition factors have labels in $K$.
These subcategories exhaust $\DD_J$.
\end{theorem}

\begin{proof}
Applying the quotient functors to endomorphisms gives the algebra maps
and their composition identities, and \eqref{eq:finite-hom-lift} proves
surjectivity. For a pair $T_L(\lambda),T_L(\mu)$, the factorization
basis of \cite[Corollary~4.44]{BS} has
$t_\lambda(\nu)t_\mu(\nu)$ maps with intermediate label $\nu$,
each factoring through $T_L(\nu)$. Those with $\nu\notin K$
are killed by $q_K^L$. They are independent, and
\eqref{eq:finite-flag-hom} counts them as the difference of the
source and target Hom dimensions. Surjectivity therefore makes
them a basis of the kernel on this Hom space. Taking every pair
of summands gives \eqref{eq:ringel-kernel}: its right-hand ideal
consists of sums of maps factoring through the discarded tiltings.

Global reconstruction of morphisms now gives
\eqref{eq:ringel-limit}, with multiplication induced by composition.
Under this isomorphism $\pi_K$ selects the $K$-coordinate.
Each $\pi_K$ is surjective: along a countable cofinal exhaustion
by finite closed subsets, one can lift recursively through the surjective
transition maps.

Surjectivity gives $\BJ/I_K\cong\BK$, so each $I_K$ has finite
codimension. The topology defined above is precisely the inverse-limit
topology: finitely many coordinate conditions are controlled by one
larger window. Thus \eqref{eq:ringel-limit} is a homeomorphism and
$\BJ$ is Hausdorff and complete. Convergence means eventual agreement
in every fixed $\BK$.

To compare with the endomorphism topology, use $q_K=e_K(-)$ to write
\[
 I_K=\{f^{\op}:f((\TJ)_\mu)=0\text{ for every }\mu\in K\}.
\]
The spaces $(\TJ)_\mu$ are finite dimensional, so this is the vanishing
of finitely many matrix coefficients. Conversely, any finite set of
coefficients involves finitely many weights. Their upper closure in
$J$ is a window, and its kernel $I_K$ forces those coefficients to
vanish. The two neighbourhood bases therefore define the same topology.
Writing $\TJ^f=\bigoplus_\mu(\TJ)_\mu^*$, the coefficient pairing is
\[
 \langle\phi\otimes v,f^{\op}\rangle=\phi(f(v)),
 \qquad \phi\otimes v\in\TJ^f\otimes_{A_J}\TJ.
\]
This is the pairing in \cite[Lemma~2.10, (2.12)--(2.15)]{BS}.
A linear functional is continuous exactly when it vanishes on some
$I_K$, so the continuous dual coalgebra is $\varinjlim_K\BK^*$.

The category $\RR_J$ is upper finite highest weight and $\TJ$
is its tilting generator. Hence \cite[Theorem~4.27]{BS} identifies
the finite-dimensional comodules over this coalgebra with the stated
lower finite highest-weight category. For finite $K\subseteq L$,
\eqref{eq:ringel-kernel} identifies $\DD_K$ with the modules over
$\mathsf B_L$ annihilated by the discarded idempotents. These are
exactly the modules with simple composition labels in $K$, giving
the asserted Serre image.

Finally, continuity of a discrete finite-dimensional module makes the
annihilator of each basis vector open. The intersection for a finite
basis contains $I_K$ for one $K$, so the action factors through $\BK$.
Conversely, an action through the discrete algebra $\BK$ is continuous.
This proves the factorization and exhaustion assertions.
\end{proof}

For $K\subseteq L$ and
$\Gamma\subseteq K$, write
$\varepsilon_\Gamma=\sum_{\lambda\in\Gamma}\varepsilon_{\lambda,L}$
and $\varepsilon_{\rm out}=\sum_{\nu\in L\setminus K}
\varepsilon_{\nu,L}$. The map induced by $\pi_{L,K}$ has kernel
\begin{equation}\label{eq:family-kernel}
 \ker\bigl(\varepsilon_\Gamma\mathsf B_L\varepsilon_\Gamma
       \longrightarrow\varepsilon_\Gamma\BK\varepsilon_\Gamma\bigr)
 =\varepsilon_\Gamma\mathsf B_L\varepsilon_{\rm out}
             \mathsf B_L\varepsilon_\Gamma.
\end{equation}
On the target side $\varepsilon_\Gamma$ denotes its image under
$\pi_{L,K}$. Taking a corner of \eqref{eq:ringel-kernel} proves the formula.

Here \emph{inflation} means keeping the vector space and using
$b\cdot v=\pi_K(b)v$; maps are unchanged. We identify $\DD_K$ with
its essential image under $i_K$. Write $\Ind(\DD_J)$ for the
category obtained by adjoining filtered direct limits.

\subsection{The Ringel functor and reciprocity}

The $A_J$-module realization gives a bimodule description of the
Ringel functor. Let
\[
 \GJ=\TJ^f\otimes_{A_J}-:
       \OO^J\longrightarrow\Ind(\DD_J).
\]
The right $A_J$-action on $\TJ^f$ is the transpose of the left action on
$\TJ$, and the commuting left $\BJ$-action is given by precomposition:
$f^{\op}\phi=\phi\circ f$. The tensor product carries the natural
comodule structure over the continuous dual coalgebra. It defines the
covariant Ringel functor, equivalently the continuous dual of
$\Hom_U(-,\TJ)$ \cite[Lemma~2.11]{BS}. Here the Hom space has
the topology of evaluation on finitely many homogeneous vectors
and linear functionals of the corresponding weight spaces of $\TJ$.
The finite counterpart is
\[
 \GK=\Hom_{\CC_K}(-,\Theta_K)^*:
             \CC_K\longrightarrow\DD_K.
\]
Both functors are right exact. On ascending standard-filtered objects,
$\GJ$ restricts to an exact equivalence onto the ascending
costandard-filtered objects of $\Ind(\DD_J)$
\cite[Theorem~4.27]{BS}.

The inflation functor $i_K$ fits into the following square, which
commutes up to natural isomorphism:
\begin{equation}\label{eq:ringel-square}
\begin{tikzcd}[column sep=large,row sep=large]
 \CC_K \arrow[r,"\GK"] \arrow[d,"j_{K,!}"'] &
 \DD_K \arrow[d,"i_K"] \\
 \OO^J \arrow[r,"\GJ"'] & \Ind(\DD_J).
\end{tikzcd}
\end{equation}
The natural isomorphism identifying the two composites follows from
$\TJ^f\otimes_{A_J}A_Je_K\cong\TJ^fe_K$ and associativity of tensor
products; it is natural in the module being tensored.

For the finite dual objects, put
$F_K=\Hom_{\CC_K}(\Theta_K,-)$. Finite Ringel duality
\cite[Definition~4.9, Theorem~4.10]{BS} gives
\[
\begin{aligned}
 P_K^\vee(\lambda)&=F_KT_K(\lambda)\cong\BK\varepsilon_{\lambda,K},
 & I_K^\vee(\lambda)&=\GK T_K(\lambda),\\
 L_K^\vee(\lambda)&=\operatorname{hd}P_K^\vee(\lambda)
                 =\operatorname{soc}I_K^\vee(\lambda),\\
 \Delta_K^\vee(\mu)&=F_K\nabla_K(\mu),
 & \nabla_K^\vee(\mu)&=\GK\Delta_K(\mu).
\end{aligned}
\]
The head and socle here are simple. Inflation identifies these finite
simples, standards and costandards with the corresponding global objects.

Write $L^\vee(\lambda)$, $\Delta^\vee(\lambda)$ and
$\nabla^\vee(\lambda)$ for the simple, standard and costandard objects
of $\DD_J$. Set
\[
 I^\vee(\lambda)=\GJ T(\lambda)
                  \in\Ind(\DD_J).
\]
It is the injective hull of $L^\vee(\lambda)$, whereas
$\GJ\Delta(\mu)=\nabla^\vee(\mu)$ belongs to $\DD_J$.

\begin{corollary}[Reciprocity]\label{thm:dual-reciprocity}
For $\lambda,\mu\in J$, the stable standard multiplicities satisfy
\begin{equation}\label{eq:global-reciprocity}
 t_\lambda(\mu)=
 [\Delta^\vee(\mu):L^\vee(\lambda)]
 = [\nabla^\vee(\mu):L^\vee(\lambda)].
\end{equation}
\end{corollary}

\begin{proof}
Choose a finite closed subset containing both labels. Finite Ringel duality
takes the standard flag of $T_K(\lambda)$ under $\GK$ to a costandard flag of the
dual injective hull. Injective reciprocity in the finite highest-weight
category identifies its $\nabla_K^\vee(\mu)$-multiplicity with
$[\Delta_K^\vee(\mu):L_K^\vee(\lambda)]$. Since inflation preserves
these composition factors, this proves the first equality. For the
second, the other finite Ringel functor,
$F_K$, takes costandards to dual standards and
$T_K(\lambda)$ to the dual projective. Projective reciprocity applies,
and self-duality of $T_K(\lambda)$ identifies its costandard
multiplicities with the same stable standard multiplicities.
\end{proof}

\subsection{The double-centralizer property}

For $\mu\in J$, define the finite dual tilting
\[
 \mathsf W_J(\mu)=\GJ Q^J(\mu)
                 \cong(\TJ)_\mu^*\in\DD_J.
\]
Put
$\mathsf W_K(\mu)=\GK Q_K(\mu)$ for $\mu\in K$ and
$\mathsf W_K=\bigoplus_{\mu\in K}\mathsf W_K(\mu)$.
The identity
$j_{K,!}Q_K(\mu)=A_Je_K\otimes_{A_K}A_K1_\mu\cong Q^J(\mu)$
and \eqref{eq:ringel-square} give
$i_K\mathsf W_K(\mu)\cong\mathsf W_J(\mu)$.

\begin{proposition}\label{prop:raw-bidual}
The Ringel functor identifies $A_J$ with the locally unital algebra
\begin{equation}\label{eq:raw-bidual}
 \bigoplus_{\mu,\nu\in J}
 \Hom_{\DD_J}\bigl(\mathsf W_J(\mu),\mathsf W_J(\nu)\bigr),
\end{equation}
where multiplication is opposite composition. It identifies
$1_\mu$ with the identity of $\mathsf W_J(\mu)$ and induces
\[
 A_K\cong\End_{\DD_K}(\mathsf W_K)^{\op}.
\]
\end{proposition}

\begin{proof}
The projectives $Q^J(\mu)$ have finite standard flags, so full
faithfulness of $\GJ$ on that class sends each map between the
$Q^J(\mu)$ to a map between their images $\mathsf W_J(\mu)$,
identifying the corresponding Hom blocks in \eqref{eq:raw-bidual}.
These identifications
preserve composition. Their direct sum therefore gives the asserted
algebra isomorphism, including its idempotents. The finite assertion
then follows from $i_K\mathsf W_K(\mu)\cong\mathsf W_J(\mu)$
and full faithfulness of $i_K$. This realizes the
double-centralizer property of
\cite[Corollary~4.30]{BS} in the weight-projective realization.
\end{proof}

An element $a\in1_\mu A_J1_\nu$ acts from $(\TJ)_\nu$ to
$(\TJ)_\mu$ and gives the transpose
$a^*:\mathsf W_J(\mu)\to\mathsf W_J(\nu)$.
Consequently the two global algebra constructions have the forms
\[
 A_J=\varinjlim_K A_K,
 \qquad
 \BJ=\varprojlim_K\BK.
\]
The first uses the corner inclusions $A_K\hookrightarrow A_L$;
the second uses the surjections $\pi_{L,K}:\mathsf B_L\to\BK$.

\subsection{Realization at finite Frobenius levels}

To realize the same dual algebra through Frobenius kernels, call
$(K,r)$ admissible if $r\geq1$ and $p^r>d(K)$. These pairs form a
directed set when ordered by inclusion of $K$ and increase of $r$.
Theorem~\ref{thm:frobenius} identifies the finite tilting endomorphism
algebras: the objects and algebras
\[
 \Theta_{r,K}=R_{r,K}\Theta_K,\qquad
 \mathsf B_{r,K}=\End_{\CC_{r,K}}(\Theta_{r,K})^{\op}
\]
satisfy
\[
 \beta_{r,K}:\BK\xrightarrow{\sim}\mathsf B_{r,K},
 \qquad f^{\op}\longmapsto(R_{r,K}f)^{\op}.
\]
The highest-weight objects of $\CC_{r,K}$ used here are transported
through this equivalence. For $K\subseteq L$ and $r\leq s$, define
\[
 \mathcal F_{(L,s),(K,r)}=H_{s,r,K}q_{s,K}^L:
 \CC_{s,L}\longrightarrow\CC_{r,K},
\]
where $H_{s,r,K}$ restricts the $G_s$-action to $G_r$, retaining the
rational torus. Restriction commutes with the support quotient, so
\[
 \mathcal F_{(L,s),(K,r)}R_{s,L}\cong R_{r,K}q_K^L.
\]
The induced End map is
$\beta_{r,K}\pi_{L,K}\beta_{s,L}^{-1}$.

\begin{corollary}\label{cor:frobenius-ringel}
There is a natural isomorphism
\begin{equation}\label{eq:frobenius-ringel-limit}
 \BJ\cong\varprojlim_{(K,r),\ p^r>d(K)}\mathsf B_{r,K}.
\end{equation}
Inflation along $\beta_{r,K}\pi_K:\BJ\to\mathsf B_{r,K}$ sends
finite-dimensional $\mathsf B_{r,K}$-modules into $\DD_J$.
For fixed $K$, their essential image consists precisely of the
objects of $i_K(\DD_K)$, for every admissible depth $r$.
\end{corollary}

\begin{proof}
Transport a compatible family $b_{r,K}$ back along $\beta_{r,K}$.
For fixed $K$, $\beta_{r,K}^{-1}(b_{r,K})$ is independent of $r$, since
any two admissible depths can be compared at a common larger depth.
Similarly, comparing two windows at a common sufficiently large depth
gives a compatible family in $\varprojlim_K\BK$. This construction
is inverse to the map induced by the global endomorphisms in
\eqref{eq:ringel-limit}. Inflation along the same maps proves the
assertion about module categories, since each $\beta_{r,K}$ is an
algebra isomorphism.
\end{proof}

Writing $\DD_{r,K}=\mathsf B_{r,K}\text{-}\mathrm{mod}_{\mathrm{fd}}$,
the finite Ringel construction itself commutes with restriction:
\begin{equation}\label{eq:frobenius-ringel-square}
\begin{tikzcd}[column sep=large,row sep=large]
 \CC_K \arrow[r,"R_{r,K}","\sim"'] \arrow[d,"\GK"'] &
 \CC_{r,K} \arrow[d,"\mathsf G_{r,K}"] \\
 \DD_K \arrow[r,"\Xi_{r,K}"',"\sim"] & \DD_{r,K}.
\end{tikzcd}
\end{equation}
Here $\mathsf G_{r,K}=\Hom_{\CC_{r,K}}(-,\Theta_{r,K})^*$,
and $\Xi_{r,K}$ gives a $\BK$-module the action
$b'\cdot v=\beta_{r,K}^{-1}(b')v$.
For $V\in\CC_K$, full faithfulness of $R_{r,K}$ identifies
$\Hom_{\CC_K}(V,\Theta_K)$ with
$\Hom_{\CC_{r,K}}(R_{r,K}V,\Theta_{r,K})$ naturally.
Taking linear duals gives the isomorphism in the square.

\subsection{Derived realization and stratification}\label{sec:derived}

The perfect-complex form of derived Ringel duality
\cite[Section~2.4]{FG} is realized by the same windows. We use finite
derived Ringel duality and Ext-full lower finite truncation
\cite[Theorems~3.59 and~4.16]{BS}.

Write $K^b$ for the bounded homotopy category and $D^b$ for the
bounded derived category, obtained by inverting quasi-isomorphisms.
The total left-derived functors $\mathbf L F$ below are computed
on projective resolutions; the exact functor $i_K$ acts termwise.

For a locally unital algebra $A$, let
$\operatorname{proj}_{\mathrm{fg}}A$ be the direct summands of finite
sums of principal projectives, and put
\[
 \Perf(A)=K^b(\operatorname{proj}_{\mathrm{fg}}A).
\]
We regard this category as a full subcategory of the derived category
of unital $A$-modules.

\begin{lemma}\label{lem:dual-derived-windows}
The functor $D^b(i_K):D^b(\DD_K)\to D^b(\DD_J)$ is fully faithful.
An object belongs to its essential image precisely when all
simple composition factors of its cohomology objects have labels in $K$.
These essential images, as $K$ varies, exhaust $D^b(\DD_J)$.
\end{lemma}

\begin{proof}
For finite $K\subseteq L$, the set $K$ is lower closed in the reversed
highest-weight order. The inclusion $\DD_K\hookrightarrow\DD_L$ is
Ext-full by \cite[Theorem~3.59(5)]{BS}: it induces isomorphisms
$\Ext^n_{\DD_K}(M,N)\cong\Ext^n_{\DD_L}(M,N)$ for every $n\geq0$.
These identify derived morphisms between objects concentrated in one
degree. Induction on the number of nonzero cohomology objects, using
the cohomological truncation triangles in each argument, then proves
full faithfulness of $D^b(\DD_K)\to D^b(\DD_L)$.

A derived morphism is represented by a roof $X^\bullet\leftarrow
Z^\bullet\to Y^\bullet$ whose left arrow is a quasi-isomorphism.
All these bounded complexes have finitely many finite-length terms.
One window contains their composition labels, as well as those in
any finite comparison of roofs witnessing equality. The finite-level
full faithfulness therefore passes to $D^b(\DD_J)$, and every bounded
complex belongs to a finite-level image.

Exactness of $i_K$ preserves cohomology and its labels. Conversely,
suppose the cohomology of a bounded complex has labels in $K$.
Its cohomological truncation triangles build it from these cohomology
objects. Full faithfulness lifts each connecting map to $D^b(\DD_K)$,
so induction places the complex in the essential image.
\end{proof}

\begin{theorem}\label{thm:derived}
Applying $\GJ$ termwise to finitely generated projective complexes
induces an equivalence
\[
 \mathbf G_J:\Perf(A_J)\xrightarrow{\sim}D^b(\DD_J).
\]
For every finite closed subset $K$, the square
\begin{equation}\label{eq:derived-square}
\begin{tikzcd}[column sep=large,row sep=large]
 \Perf(A_K) \arrow[r,"\mathbf L\GK","\sim"']
       \arrow[d,"\mathbf L j_{K,!}"'] &
 D^b(\DD_K) \arrow[d,"D^b(i_K)"] \\
 \Perf(A_J) \arrow[r,"\mathbf G_J"',"\sim"] & D^b(\DD_J)
\end{tikzcd}
\end{equation}
commutes naturally. Both vertical functors are fully faithful, and
their essential images, as $K$ varies, exhaust the global categories.
\end{theorem}

\begin{proof}
Induction sends $A_K1_\mu$ to $A_J1_\mu$ for $\mu\in K$, and the
identity $A_K=e_KA_Je_K$ identifies the Hom blocks between these
projectives before and after induction. It follows that induction is
fully faithful on their finite sums and direct summands, hence also on
bounded projective complexes. Each term of such a complex over $A_J$
is a direct summand of a finite sum of principal projectives, described
by an idempotent matrix. Choose $K$ containing all the principal
projective labels used in the finitely many terms. The idempotent
matrices and the differentials have their entries in the same corner
$e_KA_Je_K$. Thus the whole complex is induced from a bounded
projective complex over $A_K$.

The finite highest-weight algebra $A_K$ has finite global dimension,
so $\Perf(A_K)=D^b(A_K\text{-}\mathrm{mod}_{\mathrm{fd}})$.
Its highest-weight strata are copies of $\kk$, of global dimension
zero, so finite derived Ringel duality \cite[Theorem~4.16]{BS}
gives the upper equivalence in \eqref{eq:derived-square}.
The bimodule identity proving
\eqref{eq:ringel-square}, applied to projective complexes, proves
commutativity. Its bottom functor takes values in $D^b(\DD_J)$ because
$\GJ Q^J(\mu)=\mathsf W_J(\mu)$ is finite dimensional.

Any two source complexes lie in one finite window, where the finite
equivalence and the two fully faithful vertical functors identify
their morphisms. Hence $\mathbf G_J$ is fully faithful.
Lemma~\ref{lem:dual-derived-windows} also places every target object
in a finite window, and the finite equivalence there supplies a
preimage. This proves essential surjectivity.
\end{proof}

\subsubsection{The finite recovery criterion on the dual side}

Let $I_K^\vee(\lambda)$ be the injective hull of $L_K^\vee(\lambda)$
in $\DD_K$. Inside $\Ind(\DD_J)$ there are coherent inclusions
\begin{equation}\label{eq:injective-union}
 i_KI_K^\vee(\lambda)\ \subseteq\ i_LI_L^\vee(\lambda),
 \qquad
 I^\vee(\lambda)=\bigcup_{K\ni\lambda}i_KI_K^\vee(\lambda).
\end{equation}
For $\lambda\in K$, the largest subobject with labels in $K$ is
$\{v\in I^\vee(\lambda):I_Kv=0\}$. This construction is the right
adjoint to the exact inclusion on Ind categories, so it takes
$I^\vee(\lambda)$ to an injective $\BK$-module. Its socle is the
simple $L^\vee(\lambda)$ and is essential: every nonzero subobject
meets it. Hence it is the finite injective hull $i_KI_K^\vee(\lambda)$.
The inclusions follow from $I_L\subseteq I_K$ for $K\subseteq L$.
Every vector belongs to a finite-dimensional comodule, which factors
through some $\BK$; this proves the union in \eqref{eq:injective-union}.

\begin{corollary}\label{cor:perfect-recovery}
For $T=T(\lambda)$, the conditions of
Theorem~\ref{thm:recovery} are also equivalent to each of the following:
\begin{enumerate}
 \item $T$, placed in degree zero, belongs to $\Perf(A_J)$;
 \item $I^\vee(\lambda)$ has finite length;
 \item the increasing system in \eqref{eq:injective-union} stabilizes
       in a finite window.
\end{enumerate}
\end{corollary}

\begin{proof}
Each standard $\Delta(\mu)$ has a finite resolution by finitely
generated projectives. Indeed, the projective standard-flag sequence
has a kernel with a finite standard flag whose labels are strictly
higher than $\mu$. Induction terminates because $J_{\geq\mu}$ is
finite, and gives the required resolution. It follows that every
object with a finite standard flag is perfect. Conversely, a
degree-zero perfect object is finitely generated. In a bounded
projective representative, the differential onto the last positive
term is surjective and splits because that term is projective.
Removing this contractible pair and repeating eliminates the positive
acyclic tail. The remaining degree-zero term surjects onto the object.
Theorem~\ref{thm:recovery} gives the first equivalence.

The Ringel functor takes the ascending standard flag of $T$ to an
exhaustive ascending costandard flag of $I^\vee(\lambda)$ with multiplicities
$t_\lambda(\mu)$. Dual costandards have finite length, and each of
these multiplicities is finite. Thus finite standard support makes
the resulting costandard flag finite and gives finite length of the
injective. Conversely, a finite-length object cannot have infinitely
many nonzero successive sections, so finite length of the injective
forces finite standard support. Theorem~\ref{thm:recovery} proves
the second equivalence. Finally, finite length and
\eqref{eq:injective-union} imply stabilization; if the system
stabilizes, the injective belongs to a finite category and therefore
has finite length.
\end{proof}

\subsubsection{The heredity layers}

The finite Ringel algebras also express the highest-weight structure
through successive heredity ideals. Let $\lambda$ be minimal in a nonempty finite
closed subset $K$, and put $K'=K\setminus\{\lambda\}$,
$B=\BK$ and $\varepsilon=\varepsilon_{\lambda,K}$.
An ideal $I$ of $B$ is \emph{heredity} if $I^2=I$, it is projective
as a left $B$-module, and $I\rad(B)I=0$.

\begin{proposition}\label{prop:heredity-layer}
There is an exact sequence of $B$-bimodules
\begin{equation}\label{eq:heredity-sequence}
 0\longrightarrow B\varepsilon\otimes_{\kk}\varepsilon B
   \xrightarrow{\ m\ }B
   \xrightarrow{\ \pi_{K,K'}\ }\mathsf B_{K'}
   \longrightarrow0,
 \qquad m(x\otimes y)=xy.
\end{equation}
Its kernel ideal $B\varepsilon B$ is a heredity ideal. The corresponding
bounded derived recollement has factors $D^b(\DD_{K'})$ and
$D^b(\kk\text{-}\mathrm{mod}_{\mathrm{fd}})$.
\end{proposition}

\begin{proof}
Minimality gives
$T_K(\lambda)=\Delta_K(\lambda)=\nabla_K(\lambda)=L_K(\lambda)$,
so $\varepsilon B\varepsilon=\kk$. If
$m_\lambda=(\Theta_K:\Delta_K(\lambda))$, standard--costandard
orthogonality and self-duality give
$\dim B\varepsilon=\dim\varepsilon B=m_\lambda$.
The finite tilting basis identifies the kernel of $\pi_{K,K'}$ with
the maps having label $\lambda$, a space of dimension $m_\lambda^2$.
By \eqref{eq:ringel-kernel}, this kernel is $B\varepsilon B$.
Multiplication surjects onto the ideal, and its domain has the same
dimension. It is therefore an isomorphism.

The ideal is idempotent and projective on both sides. Moreover,
$\varepsilon\rad(B)\varepsilon=0$, whence
$(B\varepsilon B)\rad(B)(B\varepsilon B)=0$. Thus it is a heredity
ideal. For left modules, the central functors of the recollement are
inflation along $B\twoheadrightarrow B/B\varepsilon B=\mathsf B_{K'}$
and the functor $\varepsilon(-)$ from $B$-modules to
$\varepsilon B\varepsilon$-modules. They give the highest-weight
derived recollement \cite[Theorem~3.4, Remark~3.5(4)]{Krause}.
The stratifying-ideal conditions are also explicit in
\eqref{eq:heredity-sequence}: multiplication is an isomorphism and
higher tensor Tor over $\kk$ vanishes. Finite global dimension of
$B$ and its quotient ensures boundedness of the derived functors.
\end{proof}

Successively removing labels that are minimal in the original order
gives a finite derived stratification of every $\DD_K$, with factors
$D^b(\kk)$. Theorem~\ref{thm:derived} realizes these finite layers
in $\Perf(A_J)$.
% END ringel.tex

% BEGIN examples.tex
\section{Rank-one towers and mixed tilting algebras}
\label{sec:examples}

In rank one, the digits of the highest weight determine the local
endomorphism algebras and the depths at which their generators appear.
We first describe this tower, then extend it to mixed tilting algebras.

Throughout this section $G=SL_2$ and $p$ is any prime. We identify the
weight lattice with $\ZZ$, so that the positive root is $2$. For
$\lambda\in\ZZ$, set
\[
 J_\lambda=\{\lambda-2n:n\geq0\},\qquad
 K_{\lambda,N}=\{\lambda-2n:0\leq n\leq N\}.
\]
Here $T(\lambda)$ is taken in $\RR_{J_\lambda}$, and
$T_N(\lambda)=q_{K_{\lambda,N}}T(\lambda)$. The integer $N$ is
the root depth. Write $T^{SL_2}(m)$ for the rational tilting module
of dominant highest weight $m$.

\subsection{The endomorphism tower and its quotient maps}
Expand the integer $\lambda+1$ in $\ZZ_p$:
\[
 \lambda+1=\sum_{i\geq0}a_ip^i,\qquad 0\leq a_i<p.
\]
Put $I_\lambda=\{i:a_i\ne0\}$ and $w_i=a_ip^i$ for $i\in I_\lambda$.
For a finite subset $A\subseteq I_\lambda$, write
$d_A=\sum_{i\in A}w_i$ and $x_A=\prod_{i\in A}x_i$.
In particular $d_\varnothing=0$ and $x_\varnothing=1$. Define
\begin{equation}\label{eq:digit-model}
 R_{\lambda,N}
 =\kk[x_i:i\in I_\lambda]/(x_i^2,\ x_A:d_A>N).
\end{equation}
The square-free monomials $\{x_A:d_A\leq N\}$ form a finite basis.
For $N\leq M$, the transition $R_{\lambda,M}\to R_{\lambda,N}$
sends $x_A$ to itself if $d_A\leq N$, and to zero otherwise.

\begin{theorem}[The endomorphism tower]\label{thm:rank-one-tower}
The stable multiplicities of $T(\lambda)$ satisfy
\[
 t_\lambda(\lambda-2n)=
 \begin{cases}1&n=d_A\text{ for some finite }A\subseteq I_\lambda,\\
 0&\text{otherwise}.
 \end{cases}
\]
There are algebra isomorphisms
\[
 \End_{\CC_{K_{\lambda,N}}}(T_N(\lambda))\cong R_{\lambda,N}
\]
which identify the maps induced by
$q_{K_{\lambda,N}}^{K_{\lambda,M}}$, for $N\leq M$, with these
monomial transitions. Consequently
\[
 \End_U(T(\lambda))\cong\varprojlim_N R_{\lambda,N}.
\]
The kernel of its projection to depth $N$ consists of the formal sums
of monomials $x_A$ with $d_A>N$.
\end{theorem}

\begin{proof}
First suppose $\lambda\geq-1$ and choose $Q=p^u>\lambda+1$.
Andersen's tensor theorem \cite[Theorem~7.7(2)]{Andersen} identifies
\begin{equation}\label{eq:andersen-rank-one}
 T(\lambda)\cong
 T^{SL_2}(Q+\lambda)\otimes\Delta(-1)^{[u]},
\end{equation}
where $[u]$ denotes Frobenius twist. Since $h=2$ in type $A_1$,
the prime hypothesis $p\geq2h-2$ holds for every prime.
The signed-digit formula \cite[Proposition~5.4]{TW} gives Weyl
sections with highest weights
$Q-1+\epsilon$, where $\epsilon=\sum_i\pm a_ip^i$. The signs range
independently over the nonzero digits. For the rational Weyl
character $\chi(m)$ we have
\[
 \chi(m)=\frac{e^m-e^{-m-2}}{1-e^{-2}},\qquad
 \operatorname{ch}\Delta(-1)^{[u]}=
              \frac{e^{-Q}}{1-e^{-2Q}}.
\]
Pairing $\epsilon$ with $-\epsilon$ cancels the factor
$1-e^{-2Q}$:
\[
 \bigl(\chi(Q-1+\epsilon)+\chi(Q-1-\epsilon)\bigr)
       \frac{e^{-Q}}{1-e^{-2Q}}
 =\frac{e^{\epsilon-1}+e^{-\epsilon-1}}{1-e^{-2}}.
\]
The exponents $\epsilon-1$ are exactly $\lambda-2d_A$, so
\[
 \operatorname{ch}T(\lambda)
 =\sum_{A\subseteq I_\lambda}
       \operatorname{ch}\Delta(\lambda-2d_A).
\]
When $\lambda=-1$, the single term $\epsilon=0$ gives the same
identity directly. Since ordinary Verma characters are independent,
the displayed identity proves the multiplicity formula.

To determine multiplication, put $v=Q+\lambda+1$. For each
$i\in I_\lambda$, let $S_i$ be the stretch starting at $i$ and
ending immediately before the next nonzero digit of $v$.
These are its minimal down-admissible stretches. The loops
$L^{S_i}_{v-1}$ of \cite[Lemma~3.25]{TW} generate a commutative
square-zero algebra; their products are the loops for unions of
these stretches. Tensoring an
endomorphism with the identity of $\Delta(-1)^{[u]}$ is injective.
The source and target have dimension $2^{|I_\lambda|}$, the latter
by Proposition~\ref{prop:tilting-hom}. Hence this is an algebra
isomorphism. Choose $x_i$ to be the image of $L^{S_i}_{v-1}$.
The product $x_A$ has the down/up factorization through the rational
tilting of highest weight $v[\bigcup_{i\in A}S_i]-1
=Q+\lambda-2d_A$, by \cite[Theorem~3.2 and Lemma~3.25]{TW}.
After tensoring,
the middle module has highest weight $\lambda-2d_A$.
Its image in a window of depth $N<d_A$ is zero. Thus the high-depth
monomials are killed by the map induced by $q_{K_{\lambda,N}}$.
Surjectivity on
tilting endomorphisms and the multiplicity formula show that the
image has dimension $\#\{A:d_A\leq N\}$. The kernel is therefore
exactly the ideal in \eqref{eq:digit-model}.

Next, for $\lambda<-1$ and a fixed $N$, choose $q=p^L>\max\{N,-\lambda-1\}$
and put $\lambda'=\lambda+2q$. Fiebig's periodic equivalence
\cite[Theorem~25]{Fiebig} takes the window and its tilting to their
translates by $2q$. The first $L$ digits of $\lambda'+1$ agree with
those of $\lambda+1\in\ZZ_p$; all higher-digit monomials have depth
greater than $N$. The positive case therefore gives the required
multiplicities and quotient algebra. A finite collection of nested
windows admits one common shift, chosen using its largest depth.
Thus all its induced transition maps are
simultaneous monomial cutoffs.

Finally, choose the identifications coherently over all depths.
At the cofinal depths $N_L=p^L-1$, the model is the full square-zero
algebra on the nonzero digits with index less than $L$. Suppose a
model has been chosen through $N_L$. Compare it with a common-shift
model extending through $N_{L+1}$. Their discrepancy at $N_L$ is an
automorphism preserving every earlier cutoff ideal. The image of a
variable of depth $w$ therefore contains only monomials of depth
at least $w$, since the variable is zero at every depth below $w$.
The next model adjoins the square-zero variable $x_L$ when $a_L\ne0$.
Extend the old automorphism by fixing this variable. This is an
automorphism of the enlarged square-zero algebra, and both it and
its inverse replace monomials only by monomials of at least the same
depth. Hence every cutoff ideal is preserved. Adjusting the new
model by this extension makes it agree with all preceding choices.
Induction gives
coherent isomorphisms. The global Hom reconstruction now gives the
inverse limit and its stated kernels.
\end{proof}

\begin{corollary}\label{cor:rank-one-finite}
The global tilting $T(\lambda)$ has finite ordinary standard and
costandard flags if and only if $\lambda\geq-1$. In that case the
least depth recovering its endomorphism algebra is $N=\lambda+1$.
For $\lambda<-1$, both the kernel and cokernel of the comparison in
Theorem~\ref{thm:adjoint-comparison} are nonzero in every finite
retaining window.
\end{corollary}
\begin{proof}
For $\lambda+1\geq0$ the digits eventually vanish, and their largest
subset sum is $\lambda+1$. For $\lambda+1<0$ infinitely many digits
are nonzero. Apply Theorems~\ref{thm:rank-one-tower},
\ref{thm:recovery} and \ref{thm:adjoint-comparison}.
\end{proof}

\subsection{What the depth filtration remembers}
The dimensions of the quotients, indexed by numerical depth,
already recover the highest weight.

Put $E_{\lambda,N}=\End(T_N(\lambda))$ and
$d_\lambda(N)=\dim_\kk E_{\lambda,N}$, with $d_\lambda(-1)=0$.
The series of successive increments is
\begin{equation}\label{eq:dimension-series}
 H_\lambda(z):=\sum_{N\geq0}
       \bigl(d_\lambda(N)-d_\lambda(N-1)\bigr)z^N
   =\prod_{i\in I_\lambda}(1+z^{a_ip^i})\quad\text{in }\ZZ[[z]].
\end{equation}

\begin{theorem}[Recovery of the highest weight]\label{thm:depth-rigidity}
For fixed $p$, the dimension profile $(d_\lambda(N))_{N\geq0}$
determines $\lambda$. In particular, two numerically depth-indexed
endomorphism towers are isomorphic if and only if their labels agree.
For $L\geq1$, their prefixes through depth $p^L-1$ are isomorphic
if and only if their labels are congruent modulo $p^L$.
\end{theorem}
\begin{proof}
The profile determines $H_\lambda$, so it is enough to recover the
digits from its product expansion. If two products in
\eqref{eq:dimension-series} are equal, consider the least differing
exponent $w$ and cancel the finitely many common factors of smaller
exponent. Each is a unit in $\ZZ[[z]]$. At exponent $w$, the two remaining
products have different coefficients of $z^w$, a contradiction.
Thus the positive generator weights agree. Each of them determines
one digit: from
$w=a_ip^i$, with $1\leq a_i<p$, one recovers $i=v_p(w)$ and
$a_i=w/p^i$. All digits are determined, and
$\ZZ\hookrightarrow\ZZ_p$ is injective.

Below degree $p^L$, no digit with index at least $L$ contributes,
whereas every nonzero digit of lower index has weight at most
$(p-1)p^{L-1}<p^L$. The same cancellation recovers exactly the first
$L$ digits. Their equality is equivalently
$\lambda\equiv\mu\pmod{p^L}$. Equal digits give identical model
algebras and cutoff maps through this depth by
Theorem~\ref{thm:rank-one-tower}. Conversely an isomorphism of the
prefixes gives equality of their dimensions.
\end{proof}

For distinct labels $\lambda,\mu$, let $j=v_p(\lambda-\mu)$ and
let $a_j,b_j$ be the first different digits of $\lambda+1,\mu+1$.
The first depth at which their dimension profiles differ is
\[
 N_*=p^j\min\bigl(\{a_j,b_j\}\setminus\{0\}\bigr).
\]
Indeed, after cancelling the lower common factors, this is the least
remaining exponent. Multiplication by a series with constant term one
does not change the first degree of disagreement.

The role of numerical indexing is especially clear for negative
labels. For every $\lambda<-1$, the unfiltered topological endomorphism ring is
noncanonically isomorphic to
\[
 \varprojlim_m\kk[y_0,\ldots,y_{m-1}]/(y_0^2,\ldots,y_{m-1}^2),
\]
whose transitions set the last variable to zero. This follows by
enumerating the nonzero digits at the cofinal depths $p^L-1$.
The numerical filtration carries additional information. For example,
when $p=3$, the labels $-2$ and $-3$ have first generator depths $2$
and $1$, respectively, although their unfiltered topological rings
are isomorphic.

Finite prefixes have a corresponding limitation. For each $L$, the labels
$-2$ and $p^L-2$ give isomorphic entire prefixes through depth $p^L-1$.
The first global tilting has infinite ordinary standard support, and
the second has finite support. Thus no fixed finite prefix decides
global finiteness uniformly when the highest-weight label is unknown.

\subsection{A mixed-recovery obstruction for ordinary standards}
Recovering individual endomorphism rings need not recover maps between
different objects. The following pair of ordinary standards exhibits
this distinction.

\begin{example}\label{ex:mixed}
Let $p=3$, $J=\operatorname{down}(0)$ and $K=\{0,-2,-4,-6\}$.
Set $D_0=q_K\Delta(0)$, $D_1=q_K\Delta(-2)$ and $R=R_{1,K}$.
Then
\[
 \End_{\CC_K}(D_i)=\kk
       =\End_{\CC_{1,K}}(RD_i)\quad(i=0,1),
\]
but
\[
 \Hom_{\CC_K}(D_0,D_1)=0,\qquad
 \Hom_{\CC_{1,K}}(RD_0,RD_1)=\kk f\ne0.
\]
Consequently restriction recovers both individual End rings and fails
to be fully faithful on $\add(D_0\oplus D_1)$.
\end{example}
\begin{proof}
Use divided-power bases $z_0,z_1,z_2,z_3$ for the retained weights of
$\Delta(0)$ and $w_0,w_1,w_2$ for those of $\Delta(-2)$;
their weights are $-2i$ and $-2-2i$, respectively.
In characteristic three the degree-one operators are
\[
 \begin{array}{c|rrrr}
 &z_0&z_1&z_2&z_3\\\hline
 F&z_1&2z_2&0&0\\
 E&0&0&2z_1&z_2
 \end{array}
 \qquad
 \begin{array}{c|rrr}
 &w_0&w_1&w_2\\\hline
 F&w_1&2w_2&0\\
 E&0&w_0&0
 \end{array}.
\]
The PBW description of
$A_{1,K}$ in Section~\ref{sec:frobenius} shows that these arrows,
the torus and the degree-two divided powers determine the action.
The latter are $E^2/2$ and $F^2/2$. Operators with a peak above $0$
act as zero on both standards.

A torus-preserving endomorphism is diagonal in either basis. The
nonzero arrows in each table force all its diagonal entries to agree.
A forward map must have $z_0\mapsto0$; the $F$-arrows force
$z_1,z_2\mapsto0$, and $z_3\mapsto cw_2$ is unrestricted. Hence
$f(z_3)=w_2$, with all other values zero, spans the restricted Hom.
In the full window, $F^{(3)}z_0=z_3$ forces $c=0$.

The reverse Hom is one-dimensional before and
after restriction. Its generator is
$g(w_0)=z_1$, $g(w_1)=2z_2$, $g(w_2)=0$. The vector $z_1$ is
primitive in the full Verma module, so this is the image of a global
map. The tables give $fg=gf=0$. Thus the full two-object algebra has
one arrow $1\to0$, with vertices $0,1$ denoting $D_0,D_1$;
the restricted algebra has one arrow
in each direction with both length-two paths zero.
\end{proof}

% END examples.tex

% BEGIN mixed-algebras.tex
\subsection{Baby standards and mixed recovery}
For tiltings, both flags control the mixed Hom spaces. Continue in
type $A_1$, with positive root $2$, a fixed parity class,
$J=\mathord\downarrow b$, and
\[
 K=\{b,b-2,\ldots,a\},\qquad a=b-2N,\qquad q=p^r,\quad r\geq1.
\]
All global tiltings in this subsection belong to $\RR_J$.
For $\lambda\in J$, uniqueness identifies $T(\lambda)$ with the
tilting constructed in $\RR_{\mathord\downarrow\lambda}$.
For any $r\geq1$, let $Z_r(\eta)$
be the $G_rT$ baby Verma module induced from the positive
Borel, with highest weight $\eta$, and put
$Z_{r,K}(\eta)=q_{r,K}Z_r(\eta)$.
Write $Z^\vee_{r,K}(\eta)$ for its restricted Chevalley dual.
The restriction of an ordinary standard will have a flag by these objects.

\begin{lemma}\label{lem:baby-hw}
For every $r\geq1$, $\CC_{r,K}$ is a split highest-weight category
with standards $Z_{r,K}(\eta)$ and costandards
$Z^\vee_{r,K}(\eta)$, indexed by $\eta\in K$.
\end{lemma}
\begin{proof}
Use the PBW basis with exponents less than $p^r$ for $A_{r,K}$ from
Section~\ref{sec:frobenius}. Its elements have the form
$F^{(i)}1_\nu E^{(j)}$, with $0\leq i,j<q$, peak $\nu\in K$, and
both endpoints in $K$. Put $A=A_{r,K}$ and $e=1_b$. Only the
peak-$b$ element of degree zero belongs to $eAe$, so $eAe=\kk$.
The bases of $Ae$ and $eA$ are, respectively,
$F^{(i)}1_b$ and $1_bE^{(j)}$ with retained endpoints.
Multiplication sends their tensor basis bijectively to the peak-$b$
PBW elements. Thus
\[
 Ae\otimes_{\kk}eA\xrightarrow{\ \sim\ }AeA.
\]
As a left $A$-module, $AeA$ is therefore a direct sum of copies of
$Ae$, hence projective. Also $(AeA)^2=AeA$ and
$e\rad(A)e=\rad(eAe)=0$. Thus $AeA$ is a heredity ideal.

The quotient by $AeA$ has exactly the PBW elements with peak below $b$.
It is the corresponding algebra for support $\mathord\downarrow(b-2)$
and interval $K\setminus\{b\}$. Successive removal of maximal weights gives a
heredity chain.

Quotient $A1_\eta$ by the sum of the images of all maps
$A1_\nu\to A1_\eta$ with $\nu>\eta$. The resulting standard module
has basis $F^{(i)}1_\eta$, $0\leq i<q$, with retained endpoint.
The operators are those of $q_{r,K}Z_r(\eta)$. Chevalley duality
identifies the costandards.
\end{proof}

Let
\[
 t^K_\lambda(\eta)=(T_K(\lambda):\Delta_K(\eta)),
 \qquad
 h_{\lambda,r}(\eta)=\sum_{j\geq0}t^K_\lambda(\eta+2jq),
\]
where $t^K_\lambda$ is extended by zero off $K$.

\begin{theorem}\label{thm:mixed-gram}
The module $R_{r,K}T_K(\lambda)$ has both a $Z_{r,K}$-flag and a
$Z^\vee_{r,K}$-flag, each with multiplicities $h_{\lambda,r}$.
For $\lambda,\mu\in K$,
\begin{equation}
 \dim\Hom_{\CC_{r,K}}(R_{r,K}T_K(\lambda),R_{r,K}T_K(\mu))
   =\sum_{\eta\in K}h_{\lambda,r}(\eta)h_{\mu,r}(\eta).
      \label{eq:mixed-gram}
\end{equation}
The natural restriction map
\[
 \Hom_{\CC_K}(T_K(\lambda),T_K(\mu))
       \longrightarrow\Hom_{\CC_{r,K}}(R_{r,K}T_K(\lambda),
                                      R_{r,K}T_K(\mu))
\]
sends $f$ to $R_{r,K}(f)$ and is injective. Its cokernel has dimension
\begin{equation}
   \sum_{\eta\in K}
     \bigl(h_{\lambda,r}(\eta)h_{\mu,r}(\eta)
                -t^K_\lambda(\eta)t^K_\mu(\eta)\bigr).
      \label{eq:mixed-defect}
\end{equation}
These restriction maps preserve composition.
\end{theorem}
\begin{proof}
In the divided-power basis $z_i=F^{(i)}z_0$ of $\Delta(\nu)$,
\[
 F^{(m)}z_i=\binom{i+m}{m}z_{i+m},\qquad
 E^{(m)}z_i=\binom{\nu-i+m}{m}z_{i-m},
\]
where $z_i=0$ for $i<0$. Set $M_{-1}=0$, and let $M_j$ be the
span of $z_i$ with $0\leq i<(j+1)q$.
For $m<q$, positive divided powers lower the index. A negative
divided power crossing a $q$-block boundary has coefficient
$\binom{i+m}{m}=0$ in $\kk$, by Lucas' theorem: such an addition
has a carry among the first $r$ base-$p$ digits. Thus the $M_j$
form a $G_rT$ filtration.

The quotient $M_j/M_{j-1}$ has basis $z_{jq+i}$, $0\leq i<q$,
and highest weight $\nu-2jq$. For divided powers of degree less than
$q$, changing the upper argument of a binomial coefficient by a
multiple of $q$ leaves it unchanged modulo $p$, including when that
argument is negative. These formulas identify the quotient with
$Z_r(\nu-2jq)$. Applying the exact functor $q_{r,K}$ leaves a finite
flag. Refining an ordinary standard flag of $T_K(\lambda)$ gives the
asserted multiplicities. Exactness and compatibility with Chevalley
duality give the costandard flag as well.

Standard--costandard Hom orthogonality and higher Ext vanishing in
the category of Lemma~\ref{lem:baby-hw} give \eqref{eq:mixed-gram}.
Restriction of scalars is faithful on maps
(Proposition~\ref{prop:actual-inclusion}). Subtracting the ordinary tilting Hom
formula gives \eqref{eq:mixed-defect}; restriction preserves composition.
\end{proof}

\begin{corollary}\label{cor:sharp-family}
For a nonempty $\Gamma\subseteq K$, set
$M_\Gamma=\max_{\lambda\in\Gamma}(\lambda-a)/2$. Restriction is fully
faithful on $\add(\bigoplus_{\lambda\in\Gamma}T_K(\lambda))$ if and
only if $p^r>M_\Gamma$. This is equivalent to recovery of each of
the diagonal endomorphism rings.
\end{corollary}
\begin{proof}
If $q>M_\Gamma$, then $\eta\geq a$ and $j>0$ imply
$\eta+2jq>\lambda$ for every $\lambda\in\Gamma$. Thus no shifted
contribution occurs in $h_{\lambda,r}$, and all mixed blocks
are recovered. If $q\leq(\lambda-a)/2$, then $\lambda-2q\in K$ and
\[
 h_{\lambda,r}(\lambda-2q)>t^K_\lambda(\lambda-2q),
\]
because $t^K_\lambda(\lambda)=1$. Every other entry weakly
increases. The diagonal sum of squares in \eqref{eq:mixed-gram}
therefore increases strictly. This proves necessity and the diagonal
criterion.
\end{proof}

For the entire tilting generator the bound is $q>N$; a family near
the lower edge may be recovered earlier.

\begin{corollary}\label{cor:composite-family}
For a finite family of global rank-one tiltings with labels
$\Gamma\subseteq K$, the composite on endomorphism algebras through
$q_K$ and $R_{r,K}$ is an isomorphism exactly when
\[
 \bigcup_{\lambda\in\Gamma}S(T(\lambda))\subseteq K
 \quad\hbox{and}\quad p^r>M_\Gamma.
\]
\end{corollary}
\begin{proof}
Put $T_\Gamma=\bigoplus_{\lambda\in\Gamma}T(\lambda)$. In
\[
 \End_U(T_\Gamma)\xrightarrow{\,f\mapsto q_K(f)\,}
 \End_{\CC_K}(q_KT_\Gamma)\xrightarrow{\,g\mapsto R_{r,K}(g)\,}
 \End_{\CC_{r,K}}(R_{r,K}q_KT_\Gamma),
\]
the first map is surjective and the second injective. Their
composite is bijective precisely when both maps
are bijective. Apply Corollaries~\ref{cor:family-quotient}
and~\ref{cor:sharp-family}.
\end{proof}

\subsection{Transport of the rational tilting presentation}
Tensoring rational modules with an ordinary Verma below is
understood through the distribution algebra and its rational torus.

\begin{lemma}\label{lem:rational-transport}
Let $\Lambda$ be a finite set of integers $\lambda\geq-1$, in one
parity class, and choose $Q=p^u>\max_{\lambda\in\Lambda}(\lambda+1)$.
The functor
\[
 F_Q:\add\Bigl(\bigoplus_{\lambda\in\Lambda}T^{SL_2}(Q+\lambda)\Bigr)
       \longrightarrow
       \add\Bigl(\bigoplus_{\lambda\in\Lambda}T(\lambda)\Bigr),
 \qquad M\longmapsto M\otimes\Delta(-1)^{[u]},
\]
with $f\mapsto f\otimes1$, is an equivalence of additive categories,
after choosing the tensor isomorphisms on the displayed generators.
\end{lemma}
\begin{proof}
Andersen's tensor theorem \cite[Theorem~7.7(2), p.~1500]{Andersen}
identifies the objects. Its bound $p\geq2h-2$ holds for every prime
in type $A_1$. The rational Weyl support of
$T^{SL_2}(Q+\lambda)$ is $Q+S(T(\lambda))$, with the same
multiplicities. In particular, for each ordered pair $\lambda,\mu$,
both Hom dimensions equal
$\sum_\eta t_\lambda(\eta)t_\mu(\eta)$: the rational formula uses
the same summands with labels shifted by $Q$.
Tensoring with the identity of a nonzero vector space is injective
on maps and preserves composition. It follows that every mixed
Hom map is an isomorphism. The extension to sums and summands gives
the assertion.
\end{proof}

Let $Z_c$ be the endomorphism algebra
\[
 Z_c=\End_G\Bigl(
       \bigoplus_{\substack{0\leq m\leq c\\m\equiv c\ (2)}}T^{SL_2}(m)
              \Bigr).
\]
It is a finite full corner $eZe$ of the rational tilting algebra presented
by Tubbenhauer--Wedrich \cite[Theorem~A and Theorem~3.2]{TW}.
A vertex labelled $m$ denotes the summand $T^{SL_2}(m)$, and $e_m$
is its identity idempotent. Tubbenhauer--Wedrich write this label as
$v-1$; thus $v=m+1$. We use their downward and upward maps and relations.

\begin{theorem}\label{thm:quiver-transport}
Choose $H=p^\ell>\max(N,-a-1)$, put $\tau=2H$, and choose
$Q=p^u>b+\tau+1$. Set
\[
 c=Q+b+\tau,\qquad \alpha=Q+a+\tau,\qquad
 e_{<\alpha}=\sum_{\substack{0\leq m<\alpha\\m\equiv c\ (2)}}e_m.
\]
Then there is an algebra isomorphism
\begin{equation}\label{eq:quiver-transport}
 \BK\cong\bigl(Z_c/Z_ce_{<\alpha}Z_c\bigr)^{\op}.
\end{equation}
The vertex with rational highest-weight label $m=Q+\lambda+\tau$
corresponds to $T_K(\lambda)$: periodicity first sends $\lambda$ to
$\lambda+\tau$, and the tensor comparison uses
$T^{SL_2}(Q+\lambda+\tau)$.
For every finite chain of nested intervals with the same upper bound,
one can use common $\tau,Q$, and the transition maps are exactly the
quotient homomorphisms by the additional lower-vertex ideals.
\end{theorem}
\begin{proof}
The bound on $H$ makes all labels in $K'=K+\tau$ at least $-1$
and gives the separation required by Fiebig's periodicity
theorem \cite[Theorem~25, p.~1099]{Fiebig}. It supplies the exact equivalence
$E_\tau:\CC_K\to\CC_{K'}$, preserving standards, duality and tiltings.
The bound on $Q$ permits Lemma~\ref{lem:rational-transport} on
all retained summands. Apply to the rational tilting sum defining $Z_c$
the composite
\[
 M\longmapsto M\otimes\Delta(-1)^{[u]}
       \longmapsto q_{K'}(M\otimes\Delta(-1)^{[u]})
       \longmapsto E_\tau^{-1}q_{K'}(M\otimes\Delta(-1)^{[u]}).
\]
On morphisms it sends $f$ to
$E_\tau^{-1}q_{K'}(f\otimes1)$. After the identifications of
the retained summands, it induces
$\Phi:Z_c\to\End_{\CC_K}(\Theta_K)$.
The summands with $m<\alpha$ have highest weight
$m-Q<a+\tau$ and vanish. The retained summands become the designated
local tiltings. Lemma~\ref{lem:rational-transport} and surjectivity of
local tilting Hom prove fullness on all retained blocks.

To identify $\ker\Phi$, write $c_m(\eta)$ for the Weyl multiplicity
of weight $\eta$ in $T^{SL_2}(m)$. For a pair of rational tiltings
labelled $m,n$, the factorization basis of
\cite[Corollary~4.44]{BS} has $c_m(\eta)c_n(\eta)$ maps through
the tilting labelled $\eta$. The basis maps at labels below $\alpha$
factor through lower tiltings and span $Z_ce_{<\alpha}Z_c$.
The remaining basis maps on a retained pair are counted by the
intersection of its rational Weyl supports with $Q+K'$.
Their multiplicities are the corresponding local standard multiplicities, so
their number is the dimension of the local Hom block. Since the
map $\Phi$ is surjective, its kernel is exactly that ideal.
Taking opposites gives \eqref{eq:quiver-transport}.

For a finite chain choose the parameters for its largest interval.
Raising the lower endpoint kills exactly the additional lower
vertices. The identity $q_K^Lq_L^M=q_K^M$ identifies
these maps, including their products, with the transition maps
$\pi_{L,K}$.
\end{proof}

Thus \eqref{eq:quiver-transport} gives compatible
presentations for the finite quotients in the inverse system for $\BJ$.

\subsection{A four-weight algebra with its products}
\begin{example}\label{ex:mixed-four}
Take $p=3$ and $K=\{0,-2,-4,-6\}$. In this order the ordinary
standard multiplicity rows are
\[
 (1,1,0,0),\qquad(0,1,0,1),\qquad(0,0,1,0),\qquad(0,0,0,1).
\]
The ordinary endomorphism algebra has an isolated vertex $-4$ and
a connected component
\[
 \begin{tikzcd}[column sep=large]
 0\arrow[r,bend left=17,"d"]&
 -2\arrow[l,bend left=17,"u"]\arrow[r,bend left=17,"a"]&
 -6\arrow[l,bend left=17,"b"].
 \end{tikzcd}
\]
With $gf=g\circ f$, a complete set of relations is
\[
 du=0,\qquad ab=0,\qquad ad=0,\qquad ub=0.
\]
The loops $ud$ and $ba$ are nonzero. The connected basis consists
of the three identities, $d,u,a,b,ud,ba$; every path of length at
least three is zero. Together with the isolated identity this gives
dimension $10$. The Ringel algebra $\BK$ is its opposite.

In Theorem~\ref{thm:quiver-transport}, the choice $\tau=18,Q=27$
uses rational tiltings of highest weights $45,43,41,39$.
Their tensor factors follow from \cite[Proposition~3.4(a)]{STWZ}.
We describe the resulting maps. Put $V=\kk^2$ and
$\mathrm{St}=\operatorname{Sym}^2V$. The relevant rational factors are
$T^{SL_2}(3)=V\otimes\mathrm{St}$ and
$T^{SL_2}(4)=\operatorname{Sym}^2(\mathrm{St})$, viewed as the
symmetric summand of $\mathrm{St}\otimes\mathrm{St}$.
In fact,
\[
 \begin{aligned}
 T^{SL_2}(45)&=T^{SL_2}(3)\otimes\mathrm{St}^{[1]}
                         \otimes T^{SL_2}(4)^{[2]},\\
 T^{SL_2}(43)&=T^{SL_2}(4)\otimes T^{SL_2}(13)^{[1]},\\
 T^{SL_2}(41)&=\mathrm{St}\otimes T^{SL_2}(13)^{[1]}.
 \end{aligned}
\]
After tensoring with $\Delta(-1)^{[3]}$ and translating by $-18$,
the higher factors contribute only their highest vector at weights
in $K$. Thus the retained actions are those on
$V\otimes\mathrm{St}\otimes\mathrm{St}^{[1]}$ shifted by $-9$,
$T^{SL_2}(4)$ and $\mathrm{St}$ shifted by $-6$, and a single
highest vector $z$ of weight $-6$.

Use monomial bases $u_0,u_1$ of $V$ and $s_0,s_1,s_2$ of
$\mathrm{St}$, and write $t_i=s_i^{[1]}$. Juxtaposition in the
following lists denotes tensor product. Bases for the retained
weight spaces of $T_K(0)$ and $T_K(-2)$ are
\[
 \begin{aligned}
 (v_1,\ldots,v_7)&=(u_0s_0t_0,u_0s_0t_1,u_0s_1t_0,u_0s_2t_0,
                         u_1s_0t_0,u_1s_1t_0,u_1s_2t_0),\\
 (w_1,\ldots,w_4)&=(s_0s_0,s_0s_1+s_1s_0,s_1s_1,s_0s_2+s_2s_0).
 \end{aligned}
\]
Their weights are respectively $(0,-6,-2,-4,-2,-4,-6)$ and
$(-2,-4,-6,-6)$. On the monomial basis $y_i$ of $\operatorname{Sym}^nV$,
\[
 E^{(m)}y_i=\binom{i}{m}y_{i-m},\qquad
 F^{(m)}y_i=\binom{n-i}{m}y_{i+m}.
\]
Together with the divided-power tensor formula and Frobenius twist,
these determine the retained operators. The arrows have the
following nonzero images; all other basis images vanish:
\[
\begin{array}{c|l}
 d & v_2,v_7\mapsto w_3+w_4;\quad v_3,v_5\mapsto w_1;
                                      \quad v_4,v_6\mapsto w_2\\
 u & w_1\mapsto v_3+2v_5;\quad w_2\mapsto2v_4+v_6;\quad
       w_3\mapsto v_2+2v_7;\quad w_4\mapsto2v_2+v_7\\
 a & w_3\mapsto2z;\quad w_4\mapsto z\\
 b & z\mapsto w_3+w_4
\end{array}
\]
Substitution verifies that these maps commute with $E^{(m)},F^{(m)}$
for $m=1,2,3$, which generate the operators needed in this window.
Their products give the four zero relations, whereas
$ud(v_3)=v_3+2v_5$ and $ba(w_3)=2(w_3+w_4)$ are nonzero.
The ordinary Hom formula gives total dimension $10$, so these
maps and the four identities form the stated basis.

The quotient functor to $K'=\{0,-2,-4\}$ induces the algebra map
$\pi_{K,K'}$, whose kernel is spanned by
$e_{-6},a,b,ba$. The quotient ideal has dimension $4$ and the
quotient algebra dimension $6$. In particular the loop at $-2$
disappears while the loop at $0$ survives.

The full and level-one Hom-dimension matrices are, respectively,
\[
 \begin{pmatrix}
 2&1&0&0\\1&2&0&1\\0&0&1&0\\0&1&0&1
 \end{pmatrix},
 \qquad
 \begin{pmatrix}
 3&2&0&1\\2&2&0&1\\0&0&1&0\\1&1&0&1
 \end{pmatrix}.
\]
At level one the first multiplicity row becomes $(1,1,0,1)$,
and the other rows are unchanged; their Gram matrix gives the second
displayed matrix. An additional map
$f:R_{1,K}T_K(0)\to R_{1,K}T_K(-6)$ sends $v_2$ to $z$ and
the other $v_i$ to zero. It commutes with operators of degrees $1,2$,
but
\[
 F^{(3)}v_1=2v_2+v_7,\qquad fF^{(3)}v_1=2z,
 \qquad F^{(3)}f(v_1)=0.
\]
Thus $f$ is not a full-distribution morphism. The endomorphism
dimension increases from $10$ to $15$, and the least recovery
level is $r=2$. The diagonal at $0$ already detects this failure,
as required by Corollary~\ref{cor:sharp-family}.
\end{example}
% END mixed-algebras.tex

% BEGIN further-questions.tex
\subsection{Further questions}
For higher rank, the finite recovery criterion in
Theorem~\ref{thm:recovery} leaves the problem of deciding directly from
the highest weight and the characteristic whether a tilting has a
finite ordinary Verma flag. One would also like to compute its least
recovering finite closed subset of the fixed support from these data.
The kernels and cokernels of the comparison maps $c_K$ in
Theorem~\ref{thm:adjoint-comparison} provide a concrete approach.

The rank-one presentations in Theorem~\ref{thm:quiver-transport}
use auxiliary periodicity and tensor parameters. A further problem is
to choose generators and relations simultaneously at all truncation
depths, compatibly with every quotient map. The numerical depth
filtration should be retained: Theorem~\ref{thm:depth-rigidity} shows
that it distinguishes highest weights whose complete topological
endomorphism algebras may be isomorphic. Presentations for
endomorphism algebras of finite sums of tiltings in higher rank
would extend these explicit calculations.

Finally, which larger derived categories admit an extension of
$\Perf(A_J)\simeq D^b(\DD_J)$? The proof of
Theorem~\ref{thm:derived} represents each finite collection of bounded
complexes and morphisms in a common finite category. Extending the
equivalence to unbounded complexes calls for a replacement for this
finiteness argument. The role of the topology on $\BJ$, and compatibility
with the filtrations defined by its finite quotients, form part of
this question.
% END further-questions.tex

\clearpage
% BEGIN bibliography.tex

% END bibliography.tex

\end{document}